\documentclass[11pt]{amsart}
\usepackage{amsmath,amsthm,amsfonts,amscd,amssymb,eucal,latexsym}

\theoremstyle{definition}
\newtheorem{theorem}{Theorem}[section]
\newtheorem{corollary}[theorem]{Corollary}
\newtheorem{lemma}[theorem]{Lemma}
\newtheorem{proposition}[theorem]{Proposition}

\newtheorem{notation}[theorem]{Notation}
\newtheorem{example}[theorem]{Example}

\newtheorem{definition}[theorem]{Definition}

\newcommand{\cardn}{\text{\rm card}}
\newcommand{\vol}{\text{\rm vol}}
\newcommand{\sep}{\text{\rm sep}}
\newcommand{\spn}{\text{\rm spn}}

\newcommand{\Entp}{\text{\rm Entp}}

\newcommand{\Ad}{\text{\rm Ad}}
\newcommand{\sa}{\text{\rm sa}}
\newcommand{\clipnorm}{$_{\text{\rm c}}$Lip-norm}
\newcommand{\clipnorms}{$_{\text{\rm c}}$Lip-norms}
\newcommand{\Mdim}{\text{\rm Mdim}}
\newcommand{\Kdim}{\text{\rm Kdim}}
\newcommand{\Aut}{\text{\rm Aut}}
\newcommand{\id}{\text{\rm id}}
\newcommand{\diam}{\text{\rm diam}}
\newcommand{\spa}{\text{\rm span}}
\newcommand{\re}{\text{\rm Re}}
\newcommand{\im}{\text{\rm Im}}

\allowdisplaybreaks

\begin{document}

\title{Dimension and dynamical 
entropy for metrized $C^*$-algebras}

\author{David Kerr}
\address{Dipartimento di Matematica, Universit\`{a} di Roma 
``La Sapienza,'' P.le Aldo Moro, 2, 00185 Rome, Italy}
\email{kerr@mat.uniroma1.it}
\date{November 4, 2002}

\begin{abstract}
We introduce notions of dimension and dynamical entropy for
unital $C^*$-algebras ``metrized'' by means of {\clipnorm}s, which
are complex-scalar versions of the Lip-norms constitutive of Rieffel's 
compact quantum metric spaces. Our examples involve
the UHF algebras $M_{p^\infty}$ and noncommutative tori. In particular
we show that the entropy of a noncommutative toral automorphism with 
respect to the canonical {\clipnorm} coincides with the topological entropy 
of its commutative analogue.
\end{abstract}

\maketitle

\section{Introduction}

The idea of a noncommutative metric space was introduced by Connes
\cite{C1,C2,C3} who showed in a noncommutative-geometric context that
a Dirac operator gives rise to a metric on the 
state space of the associated $C^*$-algebra. The question of when the
topology thus obtained agrees with the weak$^*$ topology was
pursued by Rieffel \cite{MSACG,MSS}, whose line of investigation led
to the notion of a quantum metric space defined by specifying a
Lip-norm on an order-unit space \cite{GHDQMS}. This definition 
includes Lipschitz seminorms on functions over compact metric spaces 
and more generally applies to unital $C^*$-algebras, the subspaces of
self-adjoint elements of which form important examples of 
order-unit spaces. We would like to investigate here the 
structures which arise by essentially specializing and complexifying 
Lip-norms to obtain what we call ``{\clipnorm}s'' on 
unital $C^*$-algebras. We 
thereby propose a notion of dimension for {\clipnorm}ed unital 
$C^*$-algebras, along with two dynamical entropies (the second a 
measure-theoretic version of the first) which operate within the 
restricted domain of {\clipnorm}s satisfying the Leibniz rule (our 
version of noncommutative metrics).

Means for defining dimension appeared within Rieffel's work on
quantum Gromov-Hausdorff distance in \cite{GHDQMS}, where it is
pointed out that Definition 13.4 therein gives rise to possible
``quantum'' versions of Kolmogorov $\varepsilon$-entropy.
We take here a different approach which has its origin in Rieffel's
prior study of Lip-norms in \cite{MSACG,MSS}, where the total 
boundedness of the set of elements of Lip-norm and order-unit norm
no greater than $1$ was shown to be a fundamental property. This
total boundedness leads us to a definition of dimension using
approximation by linear subspaces (Section~\ref{S-dim}). 
This also makes sense for order-unit spaces, but we will 
concentrate on $C^*$-algebras, the particular geometry 
of which will play a fundamental role in our examples,
which involve the UHF algebras $M_{p^\infty}$ and noncommutative
tori (Section~\ref{S-group}). We will also show that for usual 
Lipschitz seminorms we recover the Kolmogorov dimension
(Proposition~\ref{P-Kolm}).

We also use approximation by linear subspaces to define our 
two dynamical ``product'' entropies (Section~\ref{S-entr}). The
definitions formally echo that of Voiculescu-Brown approximation
entropy \cite{Voi,Br}, but here the algebraic structure enters
in a very different way. One drawback of the Voiculescu-Brown
entropy as a noncommutative invariant is the difficulty of obtaining
nonzero lower bounds (however frequently the entropy is in fact positive) 
in systems in which the dynamical growth is not ultimately registered 
in algebraically or statistically commutative structures. We have
in mind our main examples, the noncommutative toral
automorphisms. For these we only have partial information about the
Voiculescu-Brown entropy in the nonrational case (see \cite[Sect.\ 5]{Voi} 
and the discussion in the following paragraph), and even deciding 
when the entropy is positive is a problem 
(in the rational case, i.e., when the rotation angles with respect to
pairs of canonical unitaries
are all rational, we obtain the corresponding classical value, as
follows from the upper bound established in \cite[Sect.\ 5]{Voi}
along with the fact that the corresponding commutative toral 
automorphism sits as a subsystem, so that we can apply monotonicity). 
We show that for general noncommutative toral automorphisms the product 
entropy relative to the canonical ``metric'' coincides with the topological
entropy of the corresponding toral homeomorphism
(Section~\ref{S-tori}). In analogy
with the relation between the discrete Abelian group entropy of a discrete
Abelian group automorphism and the topological entropy of its dual 
\cite{Pet}, product entropy (which is an analytic version of discrete
Abelian group entropy) may roughly be thought of  
as a ``dual'' counterpart of Voiculescu-Brown entropy, as 
illustrated by the key role played by unitaries in obtaining 
lower bounds for the former. When passing 
from the commutative to the noncommutative in an example like the
torus, the ``dual'' unitary description persists (ensuring a 
metric rigidity that facilitates computations) while
the underlying space and the transparency of the
complete order structure vanish. The shift on $M_{p^\infty}$,
on the other hand, is equally amenable to analysis from the canonical
unitary and complete order viewpoints due to the tensor product structure,
and its value can be precisely calculated for both product entropy
(Section 6) and Voiculescu-Brown entropy \cite[Prop.\ 4.7]{Voi}
(see also \cite{CS,Voi,Th,QTE,AF} for computations with respect to other 
entropies which we will discuss below).

Since we are 
not dealing with discrete entities as in the discrete Abelian group
entropy setting, product entropies will not be $C^*$-algebraic conjugacy 
invariants, but rather bi-Lipschitz $C^*$-algebraic conjugacy 
invariants (see Definition~\ref{D-Lip}). In particular, if we 
consider {\clipnorm}s arising via the ergodic action of a
compact group $G$ equipped with a length function 
(see Example~\ref{E-groupactions}), the entropies 
will be ``$G$-$C^*$-algebraic'' invariants, that is, they will be
invariant under $C^*$-algebraic conjugacies respecting the
given group actions. To put this in context, we first point out that 
there have been two basic approaches to developing $C^*$-algebraic and von 
Neumann algebraic dynamical invariants which extend classical entropies. 
While the definitions of Voiculescu \cite{Voi} and Brown \cite{Br}
are based on local approximation, the measure-theoretic 
Connes-Narnhofer-Thirring (CNT) entropy \cite{CNT} (a generalization of 
Connes-St{\o}rmer entropy \cite{CS}) and Sauvageot-Thouvenot entropy 
\cite{ST} take a physical observable viewpoint and are defined 
via the notions of an Abelian model and a stationary coupling 
with an Abelian system, respectively (see \cite{NE} for a survey). 
Because of the role played by Abelian systems in their respective
definitions, the CNT and Sauvageot-Thouvenot entropies 
(which are known to coincide on nuclear $C^*$-algebras \cite[Prop.\ 4.1]{ST}) 
function most usefully as invariants in the asymptotically
Abelian situation. For instance, their common value for noncommutative 
$2$-toral automorphisms with respect to the canonical tracial state is 
zero for a set of rotation parameters of full Lebesgue measure \cite{NT}, 
while for rational parameters the corresponding classical value is obtained
\cite{KL} and for the countable set of irrational rotation
parameters for which the system is asymptotically Abelian (at least when
restricted to a nontrivial invariant $C^*$-subalgebra generated by a pair 
of products of powers of the canonical unitaries) the value is 
positive when the associated matrix is hyperbolic 
(see \cite[Chap.\ 9]{NE}) (in this case the Voiculescu-Brown entropy
is thus also positive by \cite[Prop.\ 4.6]{Voi}). Other
entropies which are not $C^*$-algebraic or von Neumann algebraic dynamical 
invariants have been introduced in \cite{Hu,Th,AF}. The definitions of
\cite{Hu,Th} take a noncommutative open cover approach and hence are 
difficult to compute for examples like the noncommutative toral automorphisms
(see the discussion in the last section of \cite{Hu}).
In \cite{AFTA} the Alicki-Fannes entropy \cite{AF} for general
noncommutative $2$-toral automorphisms was shown to coincide with the 
corresponding classical value if the dense algebra generated by the 
canonical unitaries is taken as the special set required by the 
definition. What is particular about the product entropies is that, from 
the perspective of noncommutative geometry as exemplified in noncommutative 
tori \cite{Rie}, they provide computable quantities which reflect the 
metric rigidity but require no additional structure to 
function (i.e., they are ``metric'' dynamical invariants).

The organization of the paper is as follows.
In Section~\ref{S-lipnorms} we recall Rieffel's definition of a 
compact quantum metric space, and with this motivation then 
introduce {\clipnorms} and the
relevant maps for {\clipnorm}ed unital $C^*$-algebras,
after which we examine some examples. In Section~\ref{S-dim}
we introduce metric dimension and establish some properties,
including its coincidence with Kolmogorov dimension for
usual Lipschitz seminorms. Section~\ref{S-group} is subdivided
into two subsections in which we compute the metric dimension
for examples arising from compact group actions on the UHF
algebras $M_{p^\infty}$ and noncommutative tori, respectively. The 
two subsections of Section~\ref{S-entr} are devoted to 
introducing the two respective product entropies and
recording some properties, and in Sections~\ref{S-shift}
and \ref{S-tori} we carry out computations for the shift
on $M_{p^\infty}$ and automorphisms of noncommutative tori,
respectively.

In this paper we will be working exclusively with {\em unital} (i.e.,
``compact'') $C^*$-algebras as generally indicated.
For a unital $C^*$-algebra $A$ we denote by $1$ its unit, by
$S(A)$ its state space, 
and by $A_{\sa}$ the real vector space of self-adjoint elements
of $A$. Other general notation is introduced in Notation~\ref{N-L}, 
\ref{N-dimension}, and \ref{N-entr1}.

This work was supported by the Natural
Sciences and Engineering Research Council of Canada. I thank
Yasuyuki Kawahigashi and the operator algebra group at the
University of Tokyo for their hospitality and for
the invigorating research environment they have provided. I also
thank Hanfeng Li for pointing out an oversight in an initial draft
and the referee for making suggestions that resulted in significant 
improvements to the paper.
\medskip

\section{{\clipnorm}s on unital $C^*$-algebras}\label{S-lipnorms}

The context for our definitions of dimension and
dynamical entropy will essentially be a
specialization of Rieffel's notion of a compact quantum metric space to the 
complex-scalared domain of $C^*$-algebras. 
A compact quantum metric space is defined by specifying a Lip-norm on an 
order-unit space (see below), and this has a natural self-adjoint 
complex-scalared interpretation on a unital $C^*$-algebra in what will call 
a ``{\clipnorm}'' (Definition~\ref{D-clipnorm}). 
In fact {\clipnorm}s will make sense in more general complex-scalared
situations (e.g., operator systems), 
as will our definition of dimension (Definition~\ref{D-dim}),
but we will stick to $C^*$-algebras as these will constitute our examples 
of interest and multiplication will ultimately enter the picture when we 
come to dynamical entropy, for which the Leibniz rule will play 
an important role. 

We begin by recalling from \cite{GHDQMS} the definition of a compact quantum 
metric space. Recall that
an {\em order-unit space} is a pair $(A,e)$ consisting of a real partially 
ordered vector space $A$ with a distinguished element $e$, called 
the {\em order unit}, such that, for each $a\in A$,
\begin{enumerate}
\item[] \begin{enumerate}
\item[(1)] there exists an $r\in\mathbb{R}$ with $a\leq re$, and
\item[(2)] if $a\leq re$ for all $r\in\mathbb{R}_{>0}$ then $a\leq 0$.
\end{enumerate}
\end{enumerate}
An order-unit space is a normed vector space under the norm
$$ \| a \| = \inf \{ r\in\mathbb{R} : -re \leq a \leq re \} , $$
from which we can recover the order via the fact that $0\leq a \leq e$
if and only if $\| a \| \leq 1$ and $\| e-a \| \leq 1$. A {\em state}
on an order-unit space $(A,e)$ is a norm-bounded linear functional
on $A$ whose dual norm and value on $e$ are both $1$ (which 
automatically implies positivity). The state space of $A$ is denoted by 
$S(A)$. An important and motivating example of an order-unit space is
provided by the space of self-adjoint elements of a unital
$C^*$-algebra. In fact every order-unit space is isomorphic to some
order-unit space of self-adjoint operators on a Hilbert space (see
\cite[Appendix 2]{GHDQMS}).
Via Kadison's function representation we also see that order-unit
spaces are precisely, up to isomorphism, the dense unital subspaces
of spaces of affine functions over compact convex subsets of
topological vector spaces (see \cite[Sect.\ 1]{MSS}).

\begin{definition}[{\rm\cite[Defns. 2.1 and 2.2]{GHDQMS}}]\label{D-cqms}
Let $A$ be an order-unit space. A {\bf Lip-norm} on $A$ 
is a seminorm $L$ on $A$ such that
\begin{enumerate}
\item[] \begin{enumerate}
\item[(1)] \hspace*{1mm}for all $a\in A$ we have $L(a) = 0$ if 
and only if $a\in\mathbb{R}e$, and
\item[(2)] \hspace*{1mm}the metric $\rho_L$ defined on the state space 
$S(A)$ by
$$ \rho_L (\mu , \nu ) = \sup \{ |\mu (a) - \nu (a) | : 
a\in A \text{ and }L(a)\leq 1 \} $$
induces the weak$^*$ topology.
\end{enumerate}
\end{enumerate}
A {\bf compact quantum metric space} is a pair $(A,L)$ consisting 
of an order-unit space $A$ with a Lip-norm $L$.
\end{definition}

As mentioned above, the subspace $A_{\sa}$ of self-adjoint elements in a 
unital $C^*$-algebra $A$ forms an order-unit space, and so
we can specialize Rieffel's definition in a more or less
straightforward way to the $C^*$-algebraic context. 
We would like, however, our ``Lip-norm''
to be meaningfully defined on the $C^*$-algebra $A$ as a vector space over
the complex numbers. Such a ``Lip-norm'' should 
be invariant under taking adjoints, and thus, after introducing some 
notation, we make the following definition, which seems reasonable 
in view of Proposition~\ref{P-selfadj}. 

\begin{notation}\label{N-L}
Let $L$ be a seminorm on the unital $C^*$-algebra $A$ which is permitted
to take the value $+\infty$. We denote the 
sets $\{ a\in A : L(a) < \infty \}$ and $\{ a\in A : L(a) \leq r \}$ (for a 
given $r>0$) by $\mathcal{L}$ and $\mathcal{L}_r$ (or in some cases for 
clarity by $\mathcal{L}^A$ and $\mathcal{L}^A_r$), respectively. 
For $r>0$ we denote by $A_r$ the norm ball $\{ a\in A : \| a \| \leq r \}$. 
We write $\rho_L$ to refer to the 
semi-metric defined on the state space $S(A)$ by
$$ \rho_L (\sigma , \omega ) = \sup_{a\in\mathcal{L}_1} |\sigma (a) - 
\omega (a) | $$
for all $\sigma, \omega\in S(A)$. We write $\diam (S(A))$ to mean the
diameter of $S(A)$ with respect to the metric $\rho_L$. We say that
$L$ {\bf separates} $S(A)$ if for every pair $\sigma , \omega$ of
distinct states on $A$ there is an $a\in\mathcal{L}$ such that
$\sigma (a) \neq \omega (a)$, which is equivalent to $\rho_L$ being
a metric.
\end{notation}

\begin{definition}\label{D-clipnorm}
By a {\bf {\clipnorm}} on a unital $C^*$-algebra $A$ we mean a seminorm $L$ 
on $A$, possibly taking the value $+\infty$, such that 
\begin{enumerate}
\item[] \begin{enumerate}
\item[(i)] \hspace*{1mm}$L(a^* ) = L(a)$ for all $a\in A$ 
(adjoint invariance),
\item[(ii)] \hspace*{1mm}for all $a\in A$ we have $L(a) = 0$ if and only 
if $a\in\mathbb{C}1$ (ergodicity),
\item[(iii)] \hspace*{1mm}$L$ separates $S(A)$ and the metric $\rho_L$ 
induces the weak$^*$ topology on $S(A)$.
\end{enumerate}
\end{enumerate}
\end{definition}

\begin{proposition}\label{P-selfadj}
Let $L$ be a {\clipnorm} on a unital $C^*$-algebra $A$. Then the restriction
$L'$ of $L$ to the order-unit space $\mathcal{L} \cap A_{\sa}$ 
is a Lip-norm, and the restriction map from $S(A)$ to $S(\mathcal{L}
\cap A_{\sa})$ is a 
weak$^*$ homeomorphism which is isometric relative to the respective metrics 
$\rho_L$ and $\rho_{L'}$. Also, if $L$ is any adjoint-invariant seminorm 
on $A$, possibly taking the value $+\infty$, such that the restriction $L'$ 
to $\mathcal{L} \cap A_{\sa}$ is a Lip-norm which separates
$S(A_{\sa})\cong S(A)$, then $L$ is a {\clipnorm}, and the restriction
from $S(A)$ to $S(\mathcal{L}\cap A_{\sa})$ is a 
weak$^*$ homeomorphism which is isometric relative to the respective metrics 
$\rho_L$ and $\rho_{L'}$.
\end{proposition}

\begin{proof}
The proposition is a consequence of the fact that if $L$ is an
adjoint-invariant seminorm then for any $\sigma , \omega\in S(A)$ 
the suprema of
$$ |\sigma (a) - \omega (a) | $$
over the respective sets $\mathcal{L}_1$ and $\mathcal{L}_1 \cap
A_{\sa}$ are the same, as shown in the discussion prior to 
Definition 2.1 in \cite{GHDQMS}. The second statement of the proposition
also requires the
fact that the ergodicity of $L'$ (condition (1) of Definition~\ref{D-cqms})
implies the ergodicity of $L$, which can be seen by noting that 
if $a\in A$ and $L(a) < \infty$ then setting 
$\re (a) = (a+a^* )/2$ and $\im (a) = (a-a^* )/2i$ (the real and 
imaginary parts of $a$) we have $L' (\re (a)) = 0$ and $L' (\im (a)) = 0$ 
by adjoint invariance, so that $\re (a) , \im (a)\in\mathbb{R}1$ by
condition (1) of Definition~\ref{D-cqms}, and hence
$a = \re (a) + i \im (a) \in\mathbb{C}1$.
\end{proof}

The following proposition follows immediately from Theorem 1.8 of 
\cite{MSACG} (note that the remark following Condition 1.5 therein shows 
that this condition holds in our case). Condition (4) 
in the proposition statement will provide the basis for our definitions of 
dimension and dynamical entropy. 

\begin{proposition}\label{P-conditions}
A seminorm $L$ on a unital $C^*$-algebra $A$, possibly taking the value 
$+\infty$, is a {\clipnorm} if and only if it separates $S(A)$ and satisfies
\begin{enumerate}
\item[] \begin{enumerate}
\item[(1)] \hspace*{1mm}$L(a^* ) = L(a)$ for all $a\in A$,
\item[(2)] \hspace*{1mm}for all $a\in A$ we have $L(a) = 0$ if and only if 
$a\in\mathbb{C}1$,
\item[(3)] \hspace*{1mm}$\sup \{ | \sigma (a) - \omega (a) | : \sigma , \omega\in S(A)
\text{ and }a\in\mathcal{L}_1 \} < \infty$,
and
\item[(4)] \hspace*{1mm}the set $\mathcal{L}_1 \cap A_1$
is totally bounded in $A$ for $\| \cdot \|$.
\end{enumerate}
\end{enumerate}
\end{proposition}

When we come to dynamical entropy, {\clipnorm}s satisfying the Leibniz rule
will be of central importance, and so we also make the following definition,
which we may think of as describing one possible noncommutative analogue
of a compact metric space (cf.\ Example~\ref{E-commutative}).

\begin{definition}
We say that a {\clipnorm} $L$ on a unital $C^*$-algebra $A$ is a 
{\bf Leibniz {\clipnorm}} if it satisfies the Leibniz rule
$$ L(ab) \leq L(a)\| b \| + \| a \|L(b) $$
for all $a,b\in\mathcal{L}$.
\end{definition}

Although we do not make lower semicontinuity a general assumption for
{\clipnorm}s, it will typically hold in our examples, and has
the advantage that we can recover the restriction of $L$ to $A_{\sa}$ 
in a straightforward manner from
$\rho_L$, as shown by the following proposition, which is a consequence of
\cite[Thm.\ 4.1]{MSS} and Proposition~\ref{P-selfadj}.

\begin{proposition}\label{P-ls}
Let $L$ be a lower semicontinuous {\clipnorm} on a unital 
$C^*$-algebra $A$. Then for all $a\in A_{\sa}$ we have
$$ L(a) = \sup \{ | \sigma (a) - \omega (a) | / \rho_L (\sigma , 
\omega ) : \sigma , \omega\in S(A) \text{ and }\sigma\neq\omega \} . $$
\end{proposition}

As for metric spaces, the essential maps in our {\clipnorm} context
are ones satisfying a Lipschitz condition, which puts a uniform bound on 
the amount of ``stretching'' as formalized in the following definition,
for which we will adopt the conventional metric space
terminology (see \cite[Defn. 1.2.1]{LA}). 

\begin{definition}\label{D-Lip}
Let $A$ and $B$ be unital $C^*$-algebras with
{\clipnorm}s $L_A$ and $L_B$, respectively.
A positive unital (linear) map $\phi : A \to B$ 
is said to be {\bf Lipschitz} if there exists a $\lambda\geq 0$ such that 
$$ L_B (\phi (a)) \leq \lambda L_A (a) $$ 
for all $a\in\mathcal{L}^A$. The least such $\lambda$ is called the 
{\bf Lipschitz number} of $\phi$.
When $\phi$ is invertible and both $\phi$ and $\phi^{-1}$ are 
Lipschitz positive we say that $\phi$ is {\bf bi-Lipschitz}. If
$$ L_B (\phi (a)) = L_A (a) $$
for all $a\in A$ then we say that $\phi$ is {\bf isometric}. 
The collection of all bi-Lipschitz $^*$-automorphisms of $A$ 
will be denoted by $\Aut_L (A)$.
\end{definition}

The category of interest for dimension will be that of {\clipnorm}ed 
unital $C^*$-algebras and Lipschitz positive unital maps, with the 
bi-Lipschitz positive unital maps forming the categorical isomorphisms. 
For entropy we will incorporate the algebraic structure in the
definitions so that we will want our positive unital maps to be in fact 
$^*$-homomorphisms. We remark that, as for usual metric spaces, 
the isometric maps are too rigid to be usefully 
considered as the categorical isomorphisms, and that our 
dimension and dynamical entropies will indeed be invariant under 
general bi-Lipschitz positive unital maps and bi-Lipschitz 
$^*$-isomorphisms, respectively. We also remark that positive unital
maps are $C^*$-norm contractive \cite[Cor.\ 1]{RD}, and hence any 
bi-Lipschitz positive unital map is $C^*$-norm isometric.

The following pair of propositions capture facts pertaining to
Lipschitz maps. The first one is clear.

\begin{proposition}\label{P-composition}
Let $A, B$, and $C$ be unital $C^*$-algebras with respective
{\clipnorm}s $L_A , L_B$, and $L_C$. If $\phi : A\to B$ and $\psi : 
B\to C$ are Lipschitz positive unital maps with Lipschitz numbers
$\lambda$ and $\zeta$, respectively, then $\psi\circ\phi$ is Lipschitz
with Lipschitz number bounded by the product $\lambda\zeta$.
\end{proposition}

\begin{lemma}\label{L-inner}
If $L$ is a {\clipnorm} on a unital $C^*$-algebra $A$ and 
$a\in\mathcal{L} \cap A_{\sa}$ then denoting by $s(a)$ the infimum of
the spectrum of $a$ we have
$$ \| a - s(a) 1 \| \leq L(a) \diam (S(A)) , $$
and hence for any $\sigma , \omega\in S(A)$ we have
$$ \rho_L (\sigma , \omega ) = \sup \{ | \sigma (a) - \omega (a) | :
a\in A_{\sa} , \, L(a) \leq 1 , \text{ and } 
\| a \| \leq \diam (S(A)) \} $$
\end{lemma} 

\begin{proof}
Let $a$ be an element of $\mathcal{L} \cap A_{\sa}$ and $s(a)$ the
infimum of its spectrum. Then there are $\sigma , \omega\in S(A)$ 
such that $\sigma (a - s(a) 1) = \| a - s(a) 1 \|$ and
$\omega (a) = s(a)$. We then have 
\begin{align*}
\| a - s(a) 1 \| & =  | \sigma (a - s(a) 1) -
\omega (a - s(a) 1) | \\
&= | \sigma (a) - \omega (a) | \\
&\leq L(a) \diam (S(A)) .
\end{align*}

The second statement of the lemma follows by noting that, for
any $\sigma , \omega\in S(A)$, 
$$ \rho_L (\sigma , \omega ) = \sup \{ | \sigma (a) - \omega (a) | :
a\in A_{\sa} \text{ and }L(a) \leq 1 \} $$
(see the first sentence in the proof of 
Proposition~\ref{P-conditions}), while if $L(a) \leq 1$ then 
$\| a - s(a) 1 \| \leq\diam (S(A))$ from above,
$$ L(a - s(a) 1) = L(a) \leq 1 $$ 
by the ergodicity of $L$, and 
$$ | \sigma (a - s(a) 1) - \omega (a - s(a) 1) |
= | \sigma (a) - \omega (a) | . $$
\end{proof}

\begin{proposition}\label{P-inner}
If $L$ is a lower semicontinuous Leibniz {\clipnorm} on a unital 
$C^*$-algebra $A$ and $u\in\mathcal{L}$ is a unitary then $\Ad u$ 
is bi-Lipschitz, and the Lipschitz numbers of $\Ad u$ and its inverse 
are bounded by $2(1 + 2L(u)\diam (S(A)))$.
\end{proposition}

\begin{proof}
By the Leibniz rule and the adjoint-invariance of $L$, 
for any $a\in\mathcal{L}$ we have
$$ L(uau^* ) \leq L(u)\| a \| + L(a) + \| a \| L(u^* ) = L(a) + 
2\| a \|L(u) . $$
For any $\sigma , \omega\in S(A)$ we therefore have,
using Lemma~\ref{L-inner} for the first equality,
\begin{align*}
\lefteqn{\rho_L (\sigma\circ\Ad u , \omega\circ\Ad u)}\hspace*{0.9cm} \\
\hspace*{0.7cm}&= \sup \{ | \sigma (uau^* ) - \omega (uau^* ) | :
a\in A_{\sa} , \, L(a) \leq 1 \text{ and } \| a \| \leq \diam (S(A)) \} \\
&\leq \sup \{ | \sigma (a) - \omega (a) | : a\in A_{\sa} \text { and }
L(a) \leq 1 + 2L(u)\diam (S(A)) \} \\
&\leq (1+2L(u)\diam (S(A))) \sup \{ | \sigma (a) - \omega (a) | :
a\in A_{\sa} \text { and }L(a) \leq 1 \} \\
&= (1+2L(u)\diam (S(A))) \rho_L (\sigma , \omega ) .
\end{align*}
Since $L$ is lower semicontinuous we can thus appeal to 
Proposition~\ref{P-ls} to obtain, for any $a\in\mathcal{L} \cap A_{\sa}$,
\begin{align*}
L(uau^* ) &= \sup_{\sigma , \omega\in S(A)} 
\frac{|(\sigma\circ\Ad u)(a) - (\omega\circ\Ad u)(a)|}{\rho_L 
(\sigma , \omega )} \\
&\leq \sup_{\sigma , \omega\in S(A)} 
\frac{|(\sigma\circ\Ad u)(a) - (\omega\circ\Ad u)(a)|}{\rho_L 
(\sigma\circ\Ad u , \omega\circ\Ad u )} \\
&\hspace*{4.5cm}\ \times\sup_{\sigma , \omega\in S(A)}
\frac{\rho_L (\sigma\circ\Ad u , \omega\circ\Ad u )}{\rho_L (\sigma ,
\omega )} \\
&= L(a) (1+ 2L(u)\diam (S(A))) .
\end{align*}
Thus, for any $a\in\mathcal{L}$, setting
$\re (a) = (a+a^* )/2$ and $\im (a) = (a-a^* )/2i$ we have
\begin{align*}
L(uau^* ) &\leq L(u \re (a) u^* ) + L(u \im (a) u^* ) \\
&\leq (L(\re (a)) + L(\im (a)))(1 + 2L(u)\diam (S(A))) \\
&\leq 2 L(a) (1 + 2L(u)\diam (S(A)))
\end{align*}
using adjoint invariance.
The same argument applies to $(\Ad u)^{-1} = \Ad u^*$, and so we
obtain the result.
\end{proof}  

We conclude this section with some examples of {\clipnorm}s.

\begin{example}[commutative $C^*$-algebras]\label{E-commutative}
For a compact metric space $(X,d)$ we define the Lipschitz seminorm
$L_d$ on $C(X)$ by
$$ L_d (f) = \sup \{ |f(x) - f(y)|/ d(x,y) : x,y\in X\text{ and } 
x\neq y \} , $$
from which we can recover $d$ via the formula
$$ d(x,y) = \sup \{ |f(x) - f(y)| : f\in C(X)\text{ and } L_d (f) 
\leq 1 \} . $$
The seminorm $L_d$ is an example of a Leibniz {\clipnorm}.  
For a reference on Lipschitz seminorms and the associated Lipschitz 
algebras see \cite{LA}.
\end{example}

\begin{example}[ergodic compact group actions]\label{E-groupactions} 
For us the most important examples
of compact noncommutative metric spaces will be those which
arise from ergodic actions of compact groups, as studied by Rieffel in 
\cite{MSACG}. 
Suppose $\gamma$ is an ergodic action of a compact group $G$ on a unital 
$C^*$-algebra $A$. Let $e$ denote the identity element of $G$.
We assume that $G$ is equipped with a length
function $\ell$, that is, a continuous function 
$\ell : G \to \mathbb{R}_{\geq 0}$ such that, for all $g,h\in G$,
\begin{enumerate}
\item[] \begin{enumerate}
\item[(1)] \hspace*{1mm}$\ell (gh)\leq \ell (g) + \ell (h)$, 
\item[(2)] \hspace*{1mm}$\ell (g^{-1}) = \ell (g)$, and 
\item[(3)] \hspace*{1mm}$\ell (g) = 0$ if and only if $g=e$. 
\end{enumerate}
\end{enumerate}
The length function $\ell$ and the group action $\gamma$
combine to produce the seminorm $L$ on $A$ defined by
$$ L(a) = \sup_{g\in G\setminus \{ e \}} 
\frac{\| \gamma_g (a) - a \|}{\ell (g)} , $$
which is evidently adjoint-invariant. It is easily verified that 
$L(a) = 0$ if and only if $a\in\mathbb{C}1$.
Also, by \cite[Thm. 2.3]{MSACG}
the metric $\rho_L$ induces the weak$^*$ topology on $S(A)$, 
and the Leibniz rule is easily checked, so that
$L$ is Leibniz {\clipnorm}.
\end{example} 

\begin{example}[quotients]\label{E-quotients}
Let $A$ and $B$ be unital $C^*$-algebras and let $\phi : A \to B$ be
a surjective unital positive linear map. For instance, $\phi$ may be a
surjective unital $^*$-homomorphism or a conditional expectation, as will
be the case in our applications. Let $L$ be a {\clipnorm} on $A$. Then
$L$ induces a {\clipnorm} $L_B$ on $B$ via the prescription
$$ L_B (b) = \inf \{ L(a) : a\in A \text{ and } \phi (a) = b \} $$
for all $b\in B$. This is the analogue of restricting a metric to a
subspace. To see that $L_B$ is indeed a {\clipnorm} we observe that
the restriction of $\phi$ to $\mathcal{L}^A \cap A_{\sa}$ 
yields a surjective morphism $\mathcal{L}^A \cap A_{\sa} \to 
\mathcal{L}^B \cap B_{\sa}$ of
order-unit spaces (for surjectivity note that if $\phi (a) \in
B_{\sa}$ then $\phi \left( \frac12 (a + a^* )\right) = \phi (a)$) so that 
we may appeal to \cite[Prop.\ 3.1]{GHDQMS} to conclude that the 
restriction of $L_B$ to $\mathcal{L} \cap B_{\sa}$ is a Lip-norm,
so that $L_B$ is a {\clipnorm} by Proposition~\ref{P-selfadj} (note that
the restriction of $L_B$ to $\mathcal{L} \cap B_{\sa}$ separates
$S(B)$, since $\phi (\mathcal{L}^A ) = \mathcal{L}^B$ and the 
restriction of $L_A$ to $\mathcal{L} \cap A_{\sa}$ separates
$S(A)$ by the first part of Proposition~\ref{P-selfadj}).
\end{example}

\section{Dimension for {\clipnorm}ed unital 
$C^*$-algebras}\label{S-dim}

Let $A$ be a unital $C^*$-algebra with {\clipnorm} $L$. Recall from
Notation~\ref{N-L} our convention that $\mathcal{L}$ and $\mathcal{L}_r$
refer to the sets $\{ a\in A : L(a)<\infty \}$ and
$\{ a\in A : L(a)\leq r \}$, 
respectively. The following notation will be extensively used for the
remainder of the article.

\begin{notation}\label{N-dimension}
For a normed linear space $(X , \| \cdot \| )$
(which in our case will either be a $C^*$-algebra or a Hilbert space) 
we will denote by $\mathcal{F}(X)$ the collection of its 
finite-dimensional subspaces, and if $Y$ and $Z$ are subsets of $X$ and 
$\delta > 0$ we will write $Y\subset_\delta Z$, and say that
{\bf $Z$ approximately contains $Y$ to within $\delta$}, if for
every $y\in Y$ there is an $x\in Z$ such that $\| y - x \| < \delta$.
Using $\dim X$ to denote the vector space dimension of a subspace $X$,
for any subset $Z\subset A$ and $\delta > 0$ we set
$$ D(Z , \delta ) = \inf \{ \dim X : X\in\mathcal{F}(A)
\text{\rm{ and }} Z \subset_\delta X \} $$
(or $D(Z , \delta )=\infty$ if the set on the right is empty) and 
if $\sigma$ is a state on $A$ then we set
$$ D_\sigma (Z , \delta ) = \inf \{ \dim X : X\in\mathcal{F}
(\mathcal{H}_\sigma) \text{\rm{ and }} \pi_\sigma (Z) \xi_\sigma 
\subset_\delta X \} $$
(or $D_\sigma (Z , \delta )=\infty$ if the set on the right is empty),
with $\pi_\sigma : A \to \mathcal{B}(\mathcal{H}_\sigma )$
referring to GNS representation associated to $\sigma$, with canonical 
cyclic vector $\xi_\sigma$.
\end{notation}

\begin{proposition}\label{P-finite}
$D(\mathcal{L}_1 , \delta )$ is finite for every $\delta > 0$.
\end{proposition}

\begin{proof}
Let $a\in\mathcal{L}$, and set $\re (a) = (a+a^* )/2$ and 
$\im (a) = (a-a^* )/2i$ (the real and imaginary parts of $a$).
Let $s(\re (a))$ and $s(\im (a))$ be the infima of the spectra of
$\re (a)$ and $\im (a)$, respectively.
Using Lemma~\ref{L-inner} and the adjoint invariance of $L$ we have
\begin{align*}
\| a - (s(\re (a)) + i s(\im (a)) ) 1 \| 
&\leq \| \re (a) - s(\re (a)) 1 \| +
\| \im (a) - s(\im (a)) 1 \|  \\
&\leq L(\re (a))\diam (S(A)) + L(\im (a))\diam (S(A)) \\
&\leq 2L(a)\diam (S(A)) .
\end{align*}
Set $r = 2\,\diam (S(A))$. Since $\mathcal{L}_1 \cap A_1$ is totally
bounded by Proposition~\ref{P-conditions}, so is 
$\mathcal{L}_r \cap A_r$ by a scaling argument. 
Let $\delta > 0$. Then there is an $X\subset\mathcal{F}(A)$ 
which approximately 
contains $\mathcal{L}_r \cap A_r$ to within $\delta$, and if
$a\in\mathcal{L}_1$ then from above we have
$$ a - (s(\re (a)) + i s(\im (a))) 1 \in
\mathcal{L}_r \cap A_r $$
so that there exists an $x\in X$ with 
$$ \| a - (s(\re (a)) + i s(\im (a))) 1 - x \| < 
\delta . $$
But $(s(\re (a)) + i s(\im (a))) 1 - x \in 
\spa (X \cup \{ 1 \})$, and so we conclude that
$$ \mathcal{L}_1 \subset_\delta \spa (X \cup \{ 1 \}) . $$
Hence $D(\mathcal{L}_1 , \delta )$ is finite.
\end{proof}

In view of Proposition~\ref{P-finite} we make the following
definition. 

\begin{definition}\label{D-dim}
We define the {\bf metric dimension} of $A$ with respect to $L$ by
$$ \Mdim_L (A) = \limsup_{\delta\to 0^+} \frac{\log 
D(\mathcal{L}_1 , \delta )}{\log\delta^{-1}} . $$
\end{definition}

We may think of $D(\mathcal{L}_1 , \delta )$ 
as the {\em $\delta$-entropy} of $A$ with respect to $L$ in analogy 
with Kolmogorov $\varepsilon$-entropy \cite{KT},
and indeed when $L$ is a Lipschitz seminorm on a compact metric space
$(X,d)$ we will recover from $\Mdim_L (C(X))$ the Kolmogorov dimension 
(Proposition~\ref{P-Kolm}).

We emphasize that in using $D(\cdot , \cdot )$ in Definition~\ref{D-dim}
(and also in the definition of entropy in Section~\ref{S-entr})
we are not making any extra geometric assumptions in our
finite-dimensional approximations by linear subspaces.
For example, we are not requiring that these subspaces be images
of positive or completely positive maps which are close to the identity
on the set in question. In computing lower bounds we are thus left to rely 
on the Hilbert space geometry implicit in the $C^*$-algebraic structure,
making repeated use of Lemma~\ref{L-on} below.

\begin{proposition}
Let $A$ and $B$ be unital $C^*$-algebras with {\clipnorm}s $L_A$ and
$L_B$, respectively. Suppose $\phi : A \to B$ is a bi-Lipschitz
positive unital map. Then 
$$ \Mdim_{L_A}(A) = \Mdim_{L_B}(B) . $$
\end{proposition}

\begin{proof}
Let $\lambda > 0$ be the Lipschitz number of $\phi$. 
Then $\phi (\mathcal{L}^A_1 )
\subset\mathcal{L}^B_\lambda$, so that if $X\in\mathcal{F}(B)$
and $\mathcal{L}^B_\lambda \subset_\delta X$ then
$$ \mathcal{L}^A_1 \subset_\delta \phi^{-1}(X) $$
since $\phi$ is isometric for the $C^*$-norm (see the remark after
Definition~\ref{D-Lip}).
As a consequence
$$ D(\mathcal{L}^A_1 , \delta ) \leq 
D(\mathcal{L}^B_\lambda , \delta ) $$
and so
\begin{align*}
\Mdim_{L_A}(A) &= \limsup_{\delta\to 0^+}
\frac{\log D(\mathcal{L}^A_1 , \delta )}{\log\delta^{-1}} \\
&\leq \limsup_{\delta\to 0^+}
\frac{\log D(\mathcal{L}^B_\lambda , \delta )}{\log\delta^{-1}} \\
&= \limsup_{\delta\to 0^+}
\frac{\log D(\mathcal{L}^B_1 , 
\lambda^{-1}\delta )}{\log\lambda^{-1}\delta^{-1}}
\cdot\lim_{\delta\to 0^+}
\frac{\log\lambda^{-1}\delta^{-1}}{\log\delta^{-1}} \\
&= \limsup_{\delta\to 0^+}
\frac{\log D(\mathcal{L}^B_1 , \lambda^{-1}\delta )}{\log
\lambda^{-1}\delta^{-1}} \\
&= \Mdim_{L_B}(B) .
\end{align*}
The reverse inequality follows by a symmetric argument.
\end{proof}

The following is immediate from Definition~\ref{D-dim}.

\begin{proposition}
Let $L$ and $L'$ be {\clipnorm}s on a unital $C^*$-algebra such that
$L\leq L'$, that is, $L(a) \leq L' (a)$ for all $a\in A$. Then
$$ \Mdim_L (A) \geq \Mdim_{L'} (A) . $$
\end{proposition}

\begin{proposition}\label{P-dimquotient}
Let $A$ and $B$ be unital $C^*$-algebras, $L_A$ a {\clipnorm} on
$A$, $\phi : A \to B$ a surjective positive unital map, and
$L_B$ the {\clipnorm} on $B$ induced from $L_A$ by $\phi$. Then
$$ \Mdim_{L_B} (B) \leq\Mdim_{L_A} (A) . $$
\end{proposition}

\begin{proof}
Since $L_B$ is induced from $L_A$ (Example~\ref{E-quotients}) for any
$b\in\mathcal{L}^B_1$ there is an $a\in A$ with $\phi (a) = b$ and
$L(a) \leq 2$. Thus if $X$ is a linear subspace of 
$A$ with $\mathcal{L}^A_2 \subset_\delta X$ it follows that
$\mathcal{L}^B_1 \subset_\delta \phi (X)$. Hence
\begin{align*}
\Mdim_{L_B}(B) &= \limsup_{\delta\to 0^+}
\frac{\log D(\mathcal{L}^B_1 , \delta )}{\log\delta^{-1}} \\
&\leq \limsup_{\delta\to 0^+}
\frac{\log D(\mathcal{L}^A_2 , \delta )}{\log\delta^{-1}} \\
&= \limsup_{\delta\to 0^+}
\frac{\log D(\mathcal{L}^A_1 , 2^{-1}\delta )}{\log 2\delta^{-1}}
\cdot\lim_{\delta\to 0^+}\frac{\log 2\delta^{-1}}{\log\delta^{-1}} \\
&= \Mdim_{L_A}(A) .
\end{align*}
\end{proof}

\begin{proposition}\label{P-dimsum}
Let $A$ and $B$ be unital $C^*$-algebras with {\clipnorm}s $L_A$ and
$L_B$, respectively. Let $L$ be a {\clipnorm} on $A\oplus B$ which induces
$L_A$ and $L_B$ via the quotients onto $A$ and $B$, respectively
(see Example~\ref{E-quotients}). Then
$$ \Mdim_L (A\oplus B) = \max (\Mdim_L (A) , \Mdim_L (B)) . $$
\end{proposition}

\begin{proof}
The inequality $\Mdim_L (A\oplus B) \geq \max (\Mdim_L (A) , \Mdim_L (B))$
follows from Proposition~\ref{P-dimquotient}. To establish the reverse
inequality, let $\delta > 0$, and let $X\in\mathcal{F}(A)$ and
$Y\in\mathcal{F}(B)$ be such that $\mathcal{L}^A_1 \subset_\delta
X$ and $\mathcal{L}^B_1 \subset_\delta Y$. If $(a,b)\in
\mathcal{L}^{A\oplus B}_1 $ then $L(a)$ and $L(b)$
are no greater than $1$, and hence there
exist $x\in X$ and $y\in Y$ such that $\| x-a \| < \delta$ and
$\| y-b \| < \delta$, so that
$$ \| (x,y) - (a,b) \| < \delta . $$
Thus
$$ \mathcal{L}^{A\oplus B}_1 \subset_\delta
\spa\left( \{ (x,0) : x\in X \} \cup \{ (0,y) : y\in Y \} \right) , $$
and so we infer that
$$ D\big(\mathcal{L}^{A\oplus B}_1 , \delta \big) \leq
D\big(\mathcal{L}^A_1 , \delta \big) + D\big(\mathcal{L}^B_1 , 
\delta \big) . $$ 
For each $\delta > 0$ the sum on the right in the above display is
bounded by twice the maximum of its two summands, and so
\begin{align*}
\Mdim_L (A\oplus B) &= \limsup_{\delta\to 0^+}
\frac{\log D\big(\mathcal{L}^{A\oplus B}_1 , 
\delta \big)}{\log\delta^{-1}} \\
&\leq \max\left( \limsup_{\delta\to 0^+}
\frac{\log 2D\big(\mathcal{L}^A_1 , 
\delta \big)}{\log\delta^{-1}} ,\, \limsup_{\delta\to 0^+}
\frac{\log 2D\big(\mathcal{L}^B_1 , \delta \big)}{\log\delta^{-1}} 
\right) \\
&= \max (\Mdim_L (A) , \Mdim_L (B)) .
\end{align*}
\end{proof}

As we show in Proposition~\ref{P-Kolm} below, if $L$ is a Lipschitz
seminorm on a compact metric space $(X,d)$ then 
$\Mdim_L (C(X))$ coincides with the Kolmogorov dimension \cite{KT,SP}, 
whose definition we recall. 
Let $(X,d)$ be a compact metric space. A set $E\subset X$ is said to be 
{\em $\delta$-separated} if for any distinct 
$x,y\in E$ we have $d(x,y) > \delta$, while a set $F\in X$ is said to 
be {\em $\delta$-spanning} if for any $x\in X$ 
there is a $y\in F$ such that $d(x,y) \leq \delta$. We denote by 
$\sep (\delta ,d)$ the largest cardinality of an $\delta$-separated
set and by $\spn (\delta ,d)$ the smallest cardinality of a
$\delta$-spanning set. We furthermore denote by
$N(\delta ,d)$ the minimal cardinality of a cover of $X$
by $\delta$-balls. The {\em Kolmogorov dimension} of $(X,d)$, which
we will denote by $\Kdim_d (X)$, is
the common value of the three expressions
$$ \limsup_{\delta\to 0^+} \frac{\log\sep (\delta ,d)}{\log\delta^{-1}} ,
\quad \limsup_{\delta\to 0^+} \frac{\log\spn (\delta ,d)}{\log\delta^{-1}} ,
\quad \limsup_{\delta\to 0^+} 
\frac{\log N(\delta ,d)}{\log\delta^{-1}} . $$
This also goes by other names in the literature, such as box dimension and 
limit capacity (see \cite[Chap.\ 2]{Pes}).  

We will need the following lemma from \cite{Voi},
which will also be of use later on.

\begin{lemma}[{\rm\cite[Lemma 7.8]{Voi}}]\label{L-on}
If $B$ is an orthonormal set of vectors in a Hilbert space $\mathcal{H}$
and $\delta > 0$ then
$$ \inf \{ \dim X : X\in\mathcal{F}(\mathcal{H})\text{ and }X
\subset_\delta B \} \geq (1-\delta^2 )\cardn (B) . $$
\end{lemma}

\begin{proposition}\label{P-Kolm}
Let $(X,d)$ be a compact metric space, and let $L$ be the associated 
Lipschitz seminorm on $C(X)$, that is,
$$ L(f) = \sup\left\{ |f(x) - f(y)| / d(x,y) : x,y\in X
\text{ and }x\neq y \right\} $$
for all $f\in C(X)$. Then
$$ \Mdim_L (C(X)) = \Kdim_d (X) . $$
\end{proposition}

\begin{proof}
Let $\delta > 0$ and let $\mathcal{U} = \{ \mathcal{B}(x_1 , \delta ),
\dots , \mathcal{B}(x_r , \delta ) \}$ be a cover of $X$ by 
$\delta$-balls. Let $\Omega = \{ f_1 , \dots , f_r \}$ be a partition of 
unity 
subordinate to $\mathcal{U}$. If $f\in\mathcal{L}_1$ and $x$ and $y$
are points of $X$ contained in the same member of $\mathcal{U}$, then
$$ |f(x) - f(y)| < 2\delta . $$
Thus for any $x\in X$ we have
\begin{align*}
\Bigg|\, f(x) - \sum_{1\leq i \leq r} f(x_i ) f_i (x) \,\Bigg| 
&\leq \sum_{1\leq i \leq r} | f(x) - f(x_i ) | \, f_i (x) \\
&\leq \sum_{\{ i : x\in\mathcal{B}(x_i , \delta )\}} | f(x) - f(x_i ) | 
\, f_i (x) \\
&< 2\delta .
\end{align*}
Thus $\mathcal{L}_1 \subset_{2\delta} \spa (\Omega )$, and since
$\dim (\spa (\Omega )) = \text{card}(\mathcal{U})$ we conclude that
$$ \Mdim_L (C(X)) = \limsup_{\delta\to 0^+}\frac{\log D(\mathcal{L}_1 , 
2\delta )}{\log\delta^{-1}} \leq \Kdim_d (X)) . $$

To establish the reverse inequality, let $\delta > 0$ and let
$E = \{ x_1 , \dots , x_r \}$ be a $\delta$-separated set of maximal
cardinality. The idea will be to consider the probability measure
$\mu$ uniformly supported on $E$ and to construct unitaries in $C(X)$ 
with sufficiently small Lipschitz seminorm which, when viewed as elements
of $L^2 (X, \mu )$, form an orthonormal basis, so that we can appeal
to Lemma~\ref{L-on}. For each $j=1, \dots , r$ define the function $f_j$ 
by $$ f_j (x) = \max (0, 1-\delta^{-1} d(x,x_j )) $$
for all $x\in X$, and observe that $L(f_j ) = \delta^{-1}$. For each
$k=1, \dots , r$ define the function $g_k$ by
$$ g_k = \sum_{j=1}^n \left[ jkr^{-1} \right] f_j , $$
where $\left[ \,\cdot\, \right]$ means take the fractional part. We then
have $L(g_k ) \leq \delta^{-1}$, as can be seen by alternatively expressing
$g_k$ as the join of the functions $\left[ jkr^{-1} \right] f_j$ (note that
the supports of the $f_j$'s are pairwise disjoint) and applying the
inequality $L(f \vee g)\leq \max (L(f) , L(g))$ relating $L$ to the lattice
structure of real-valued functions on $X$. For each $k=1 , \dots , r$
set
$$ u_k = e^{2\pi i g_k} . $$
Repeated application of the Leibniz rule yields, for each $n\geq 1$,
\begin{align*}
L\Bigg( \sum_{j=0}^n \frac{(2\pi i g_k )^j}{j!} \Bigg) \leq
\sum_{j=0}^n \frac{(2\pi )^j}{j!} L(g_k^j ) 
&\leq \sum_{j=0}^n \frac{(2\pi )^j}{j!} j L(g_k ) \\
&= 2\pi \Bigg( \sum_{j=0}^{n-1}\frac{(2\pi )^j}{j!} \Bigg) L(g_k ) \\
&\leq 2\pi e^{2\pi}L(g_k ) , 
\end{align*}
and thus, since the sequence $\left\{ \sum_{j=0}^n 
\frac{(2\pi i g_k )^j}{j!} \right\}_{n\in\mathbb{N}}$ converges uniformly
to $u_k$, we can appeal to the lower semicontinuity of $L$ to obtain the
estimate
$$ L(u_k ) \leq 2\pi e^{2\pi}L(g_k ) \leq 2\pi e^{2\pi}\delta^{-1} . $$
Setting $U(\delta ) = \{ u_k : k=1, \dots , r \}$ and $C=2\pi e^{2\pi}$, 
we thus have that the set $\{ C^{-1}u : u\in U(\delta ) \}$, which we will
simply denote by $C^{-1}U(\delta )$, lies in $\mathcal{L}_1$
if $\delta\leq C$.

Next, let $\mu$ be the probability measure uniformly supported on $E$
and let $\pi_\mu : C(X) \to \mathcal{B}(L^2 (X,\mu ))$ be the associated
GNS representation, with canonical cyclic vector $\xi_\mu$. Then, for each
$k=1, \dots , r$, $\pi_\mu (u_k ) \xi_\mu$ is the unit vector
$$ (1, e^{2\pi i kr^{-1}}, (e^{2\pi i kr^{-1}})^2 , \dots ,
(e^{2\pi i kr^{-1}})^{r-1}) $$
under the obvious identification of $L^2 (X,\mu )$ with $\mathbb{C}^r$
which respects the order of the indexing of the points $x_1 , \dots , x_r$.
Hence we see that the set $\{ \pi_\mu (u) \xi_\mu : u\in U(\delta ) \}$
forms an orthonormal basis for $L^2 (X,\mu )$, and so by Lemma~\ref{L-on}
we have
$$ D_\mu (U(\delta ) , 2^{-1}) \geq (1-2^{-2})\,\text{card}
(U(\delta )) = \textstyle{\frac34}\text{card}(E) =  \textstyle{\frac34}
\sep (\delta , d) , $$
(for the meaning of $D_\mu (\cdot , \cdot )$ see Notation~\ref{N-dimension}).

Carrying out the above construction for each $\delta > 0$, we then have
\begin{align*}
\Mdim_L (C(X)) &= \limsup_{\delta\to 0^+}
\frac{\log D(\mathcal{L}_1 , 2^{-1}C^{-1}\delta )}{\log 2C
\delta^{-1}} \\
&= \limsup_{\delta\to 0^+}\frac{\log D(\mathcal{L}_1 , 
2^{-1}C^{-1}\delta )}{\log\delta^{-1}} \\
&\geq \limsup_{\delta\to 0^+}
\frac{\log D(C^{-1}\delta U(\delta ) , 2^{-1}C^{-1}
\delta )}{\log\delta^{-1}} \\
&= \limsup_{\delta\to 0^+}\frac{\log D(U(\delta ) , 
2^{-1})}{\log\delta^{-1}} \\
&\geq \limsup_{\delta\to 0^+}\frac{\log D_\mu (U(\delta ) , 
2^{-1})}{\log\delta^{-1}} \\
&\geq \limsup_{\delta\to 0^+}\frac{\log\textstyle\frac34 \,
\sep (\delta , d)}{\log\delta^{-1}} \\
&= \Kdim_d (X) .
\end{align*}
\end{proof}

\section{Group actions and dimension}\label{S-group}

Here we compute the dimension for some examples in which the {\clipnorm} 
is defined by means of an ergodic compact group action.

\subsection{The UHF algebra $M_{p^\infty}$}\label{SS-UHF}
We consider here the infinite tensor product $M_p^{\otimes\mathbb{Z}}$
(usually denoted $M_{p^\infty}$) of $p\times p$ matrix algebras 
$M_p$ over $\mathbb{C}$ with the infinite product of Weyl actions.
As shown in \cite{OPT} there is a unique ergodic action
of $G = \mathbb{Z}_p \times\mathbb{Z}_p$ on a simple $C^*$-algebra up
to conjugacy, namely the Weyl action on $M_p$, defined as follows.
Let $\rho$ be the $p$th root of unity $e^{2\pi i/p}$, and consider
the unitary $u=\text{diag}(1, \rho , \rho^2 , \dots , \rho^{p-1})$
along with the unitary $v$ which has $1$'s on the superdiagonal and in the
bottom left-hand entry and $0$'s elsewhere. Then we have
$$ vu = \rho uv , $$
and $u$ and $v$ generate $M_p$ $C^*$-algebraically. The Weyl action
$\gamma : G \to \Aut (M_p )$ is 
given by the following specification on the generators $u$ and $v$:
\begin{align*}
\gamma_{(r,s)}(u) &= \rho^r u , \\
\gamma_{(r,s)}(v) &= \rho^s v .
\end{align*}
We may then consider the infinite product action 
$\gamma^{\otimes\mathbb{Z}}$ of the product group 
$G^{\mathbb{Z}}$ on $M_p^{\otimes\mathbb{Z}}$.

Consider the metric on $G$ obtained by viewing $G$ 
as a subgroup of $\mathbb{R}^2 / \mathbb{Z}^2$ with the metric induced 
from the Euclidean metric on $\mathbb{R}^2$, and let $\ell$ be the 
length function on $G$ defined by taking the distance to $0$. 
Given $0 < \lambda < 1$ we define the length function $\ell_\lambda$ on
$G^{\mathbb{Z}}$ by
$$ \ell_\lambda ((g_j , h_j )_{j\in\mathbb{Z}}) = \sum_{j\in\mathbb{Z}}
\lambda^{|j|} \ell ((g_j , h_j )) . $$
We could also define length functions on $G^{\mathbb{Z}}$ by using suitable
choices of weightings of $\ell$ on the factors other than the above geometric
ones (and in many cases compute the metric dimension as in 
Proposition~\ref{P-Mpdim} below), 
but for simplicity we will restrict our attention to length functions of
the form $\ell_\lambda$.
Let $L$ be the {\clipnorm} on $M_p^{\otimes\mathbb{Z}}$ arising from 
the action $\gamma^{\otimes\mathbb{Z}}$ and the length function 
$\ell_\lambda$. 

\begin{proposition}\label{P-Mpdim}
We have
$$ \Mdim_L \big( M_p^{\otimes\mathbb{Z}}\big) = 
\frac{4\log p}{\log\lambda^{-1}} . $$
\end{proposition}
   
\begin{proof}
For each $n$ consider the conditional expectation 
$E_n$ of $M_p^{\otimes\mathbb{Z}}$ onto the subalgebra 
$M_p^{\otimes [-n,n]}$ given by
$$ E_n (a) = \int_{G^{\mathbb{Z}\setminus [-n,n]}} 
\gamma_{g}^{\otimes\mathbb{Z}}(a)\, dg , $$
where $dg$ is normalized Haar measure on $G^{\mathbb{Z}}$ and
$G^{\mathbb{Z}\setminus [-n,n]}$ is the subgroup of $G^{\mathbb{Z}}$ of
elements which are the identity at the coordinates in the interval
$[-n , n]$. Then for each $a\in\mathcal{L}$ we have
\begin{align*}
\| E_n (a) - a \| &= \bigg\| \int_{G^{\mathbb{Z}\setminus [-n,n]}} 
\big( \gamma_{g}^{\otimes\mathbb{Z}}(a) - a\big) \, dg \bigg\| \\ 
&\leq \int_{G^{\mathbb{Z}\setminus [-n,n]}}
\big\| \gamma_{g}^{\otimes\mathbb{Z}}(a) - a \big\| \, dg \\
&\leq \int_{G^{\mathbb{Z}\setminus [-n,n]}} L(a) \ell_\lambda (g)
\, dg \\
&\leq L(a) \frac{2\lambda^{n+1}}{1-\lambda} .
\end{align*}
Let $\delta > 0$. If $\delta$ is sufficiently small
there is an $n\in\mathbb{N}$ such that
$$ 2\lambda^{n+1}(1-\lambda )^{-1} \leq \delta \leq 
2\lambda^n (1-\lambda )^{-1} $$
Then, in view of the above estimate on $\| E_n (a) - a \|$ when 
$a\in\mathcal{L}_1$, we have that $\mathcal{L}_1$ 
is approximately contained in
$M_p^{\otimes [-n,n]}$ to within $\delta$. Since $M_p^{\otimes [-n,n]}$
has linear dimension $p^{2(2n+1)}$ it follows that
\begin{align*} 
\frac{\log D(\mathcal{L}_1 , \delta )}{\log\delta^{-1}} &\leq 
\frac{\log D(\mathcal{L}_1 , 2\lambda^{n+1} 
(1-\lambda )^{-1})}{\log (2(1-\lambda )\lambda^{-n})} \\
&\leq \frac{(4n+2)\log p}{\log (2(1-\lambda )\lambda^{-n})} 
\end{align*}
and so
\begin{align*}
\Mdim_L \big( M_p^{\otimes\mathbb{Z}}\big) &= \limsup_{n\to\infty}
\frac{\log D(\mathcal{L}_1 , \delta )}{\log\delta^{-1}} \\
&\leq \lim_{n\to\infty}\frac{(4n+2)\log p}{\log (2(1-\lambda )
\lambda^{-n})}\\
&= \frac{4\log p}{\log\lambda^{-1}} .
\end{align*}

To prove the reverse inequality, consider for each $n\in\mathbb{N}$
the subset 
\begin{gather*}
U_n = \big\{ u^{i_{-n}}v^{j_{-n}} \otimes 
u^{i_{-n+1}}v^{j_{-n+1}} \otimes\cdots\otimes u^{i_n}v^{j_n} :\hspace*{3cm}\\
\hspace*{5cm}0\leq i_k ,j_k \leq p-1 \text{ for } k=-n , \dots , n \big\}
\end{gather*} 
of $M_p^{\otimes [-n,n]}$ (i.e., all elementary tensors in
$M_p^{\otimes [-n,n]}$ whose components are Weyl unitaries in the respective
copies of $M_p$). It is easily checked that the {\clipnorm} of any element
in $U_n$ is bounded by $2\big( 1 + 2\sum_{k=1}^n \lambda^k \big) 
\leq (4n+2)\lambda^n$. Now the product of any two distinct products
of powers of Weyl generators in $M_p$ is zero under evaluation at the
unique tracial state $\tau$ on $M_p^{\otimes\mathbb{Z}}$, as can be seen
from the commutation relation between $u$ and $v$. Thus, since 
$\tau$ is a tensor product of traces in its restriction to 
$M_p^{\otimes [-n,n]}$, the
product of any two distinct elements of $\Omega_n$ is zero under
evaluation by $\tau$. This implies that $\pi_\tau (U_n )\xi_\tau$
is an orthonormal set in the GNS representation Hilbert space 
associated to $\tau$ with canonical cyclic vector $\xi_\tau$, and
so by Lemma~\ref{L-on} we have
$$ D_\tau (U_n , 2^{-1}) \geq (1-2^{-1}) \cardn 
(\pi_\tau (U_n )\xi_\tau ) = \textstyle{\frac34} p^{2(2n+1)} . $$ 
Thus setting 
$$ W_n = \{ (4n+2)^{-1}\lambda^n w : w\in U_n \} $$
(which is contained in $\mathcal{L}_1$) 
we have 
\begin{align*}
D(W_n , (4n+2)^{-1}\lambda^{-n}2^{-1}) &\geq
D_\tau (W_n , (4n+2)^{-1}\lambda^{-n}2^{-1}) \\
&\geq D_\tau (U_n , 2^{-1}) \\
&\geq \textstyle{\frac34} p^{2(2n+1)}
\end{align*}
and so
\begin{align*}
\Mdim_L \big( M_p^{\otimes\mathbb{Z}}\big) &\geq \limsup_{n\to\infty} 
\frac{\log (D(W_n , (4n+2)^{-1}\lambda^{-n}2^{-1})}{\log
(4n+2)^{-1}\lambda^{-n}2^{-1})} \\
&\geq \limsup_{n\to\infty} 
\frac{\log\textstyle{\frac34} + (4n+2)\log p}{\log ((4n+2)^{-1} 2^{-1}) +
n\log \lambda^{-1}} \\
&= \frac{4\log p}{\log\lambda^{-1}} , 
\end{align*}
completing the proof.
\end{proof}

Because we have used the
canonical unitary desciption of $M_p^{\otimes\mathbb{Z}}$ in an essential
way, we cannot expect to be able to carry out a computation for much more
general types of 
tensor products by extending the arguments of this subsection, although 
such a computation would be possible, for example, for tensor products 
of noncommutative tori, in which case we could incorporate the methods of 
the next subsection.

\subsection{Noncommutative tori}\label{SS-noncommtorusdim}
Let $\rho : \mathbb{Z}^p \times \mathbb{Z}^p \to \mathbb{T}$ be an
antisymmetric bicharacter and for $1\leq i,j \leq k$ set
$$ \rho_{ij} = \rho (e_i , e_j ) $$
where $\{ e_1 , \dots , e_p \}$ is the standard basis for $\mathbb{Z}^p$.
The universal $C^*$-algebra $A_\rho$ generated by unitaries
$u_1 , \dots , u_p$ satisfying
$$ u_j u_i = \rho_{ij} u_i u_j $$
is referred to as a {\em noncommutative $p$-torus}. Slawny showed
in \cite{Sla} that $A_\rho$ is simple if and only if $\rho$ is
nondegenerate (meaning that $\rho (g,h) =1$ for all $h\in\mathbb{Z}^p$
implies that $g=0$), and these two conditions are furthermore
equivalent to the existence of a unique tracial state on $A_\rho$
(see \cite{Ell}).

Let $A_\rho$ be a noncommutative $p$-torus with generators
$u_1 , \dots , u_p$. There is an ergodic action $\gamma : \mathbb{T}^p
\cong (\mathbb{R} / \mathbb{Z})^p \to \Aut (A_\rho )$ determined by
$$ \gamma_{(t_1 , \dots , t_p )}(u_j ) = e^{2\pi it_j} u_j $$
(see \cite{OPT}). We will consider the {\clipnorm} $L$ arising from the
action $\gamma$
as in Example~\ref{E-groupactions}, with the length function given
by taking the distance to $0$ with respect to the metric induced from
the Euclidean metric on $\mathbb{R}^p$ scaled by $2\pi$ (scaling will
not affect the value of $\Mdim_L (A_\rho )$ but our choice of length
function ensures for convenience that $L(u_j ) = 1$ for each $j=1, \dots ,p$).
We denote by $\tau$ the tracial state defined by
$$ \tau (a) = \int_{\mathbb{T}^p} 
\gamma_{(t_1 , \dots , t_p )}(a) \, d(t_1 , \dots , t_p ) $$
for all $a\in A_\rho$, where $d(t_1 , \dots , t_p )$ is normalized Haar
measure on $\mathbb{T}^p \cong (\mathbb{R} / \mathbb{Z})^p$.

For $(n_1 , \dots , n_p )\in\mathbb{N}^p$ let
$R(n_1 , \dots , n_p )$ denote the set of points $(k_1 , \dots ,
k_p )$ in $\mathbb{Z}^p$ such that $|k_i | \leq n_i$ for
$i=1, \dots , p$. For each $a\in A_\rho$, we 
define for each $(n_1 , \dots , n_p )\in\mathbb{N}^p$ the partial 
Fourier sum
$$ s_{(n_1 , \dots , n_p )} (a) = \sum_{(k_1 , \dots , k_p )\in
R(n_1 , \dots , n_p )} \tau (au_p^{-k_p}\cdots
u_{1}^{-k_1}) u_{1}^{k_1}\cdots u_p^{k_p} $$
and for each $n\in\mathbb{N}$ the Ces\`{a}ro mean
$$ \sigma_n (a) = \left( \sum_{(n_1 , \dots , n_p )\in R(n, n, \dots , n)}
s_{(n_1 , \dots , n_p )}(a) \right) \big/ (n+1)^p . $$
Weaver showed in \cite[Thm.\ 22]{LADVNA} for the case $p=2$ 
that $\sigma_n (a)\to a$ in norm for all 
$a\in\mathcal{L}$. To compute $\Mdim_L (A_\rho )$ we will need a 
handle on the rate of
this convergence, and so we have in Lemma~\ref{L-convrate} below an 
extension to the noncommutative case of a standard result in classical 
Fourier analysis (see for example \cite{Kat}). To make the required 
estimate we will use the expression for 
$\sigma_n (a) - a$ given by the following lemma, which can be proved in
the same way as its specialization to the case $p=2$, which 
appears in a more general form in \cite{LANT} as Lemma 3.1 and is established 
in the course of the proof of \cite[Thm.\ 22]{LADVNA}.

Recall the classical Fej\'{e}r kernel $K_n$ defined by 
$$ K_n (t) = \sum_{k=-n}^n \left( 1 - \frac{|k|}{n+1} \right) e^{2\pi ikt}
= \frac{1}{n+1} \left( \frac{\sin ((n+1)t/2 )}{\sin (t/2 )} \right)^2 . $$

\begin{lemma}\label{L-integrals}
If $a\in A_\rho$ then for all $n\in\mathbb{N}$ we have
$$ \sigma_n (a) = \int_{\mathbb{T}^p}
\gamma_{(t_1 , \dots , t_p )}(a) K_n (t_1 ) \cdots 
K_n (t_p )\, d(t_1 , \dots , t_p ) $$
and
\begin{gather*}
a - \sigma_n (a) = \sum_{k=1}^p \int_{\mathbb{T}^{k-1}} 
\gamma_{(t_1 , \dots , t_{k-1}, 0 , \dots , 0)} \left( \int_{\mathbb{T}}
(a - \gamma_{r_k (t_k)}(a)) K_n (t_k ) dt_k \right) \hspace*{1cm} \\
\hspace*{5cm} \times K_n (t_1 ) \cdots K_n (t_{k-1})
\, d(t_1 , \dots , t_{k-1}) , 
\end{gather*}
with the integrals taken in the Riemann 
sense and $r_k (t)$ denoting the $p$-tuple which is $t$ at the $k$th 
coordinate and $0$ elsewhere.
\end{lemma}

Notice that the right-hand expression in the second display in 
the statement of the above lemma is a telescoping sum, so that the 
second display is an immediate consequence of the first display 
in view of the fact that the integral of the Fej\'{e}r kernel over 
$\mathbb{T}$ is $1$. Note also that the first display
shows that $\| \sigma_n (a) \| \leq \| a \|$ for all $n\in\mathbb{N}$ and
$a\in A_\rho$, a fact which will be of use in the proof of 
Proposition~\ref{P-upperbound}.

\begin{lemma}\label{L-convrate}
If $a\in\mathcal{L}^{A_\rho}$ then there is a $C>0$ such that 
$$ \| a - \sigma_n (a) \| < L(a) C \frac{\log n}{n} $$ 
for all $n\in\mathbb{N}$.
\end{lemma}

\begin{proof}
It suffices to show that each of the summands on the right-hand
side of the second display of Lemma~\ref{L-integrals} is bounded
by $M n^{-1}\log n$ for some $M>0$ and all $n\in\mathbb{N}$.
We thus observe that if $1\leq k \leq p$ then, with $r_k (t)$ denoting the 
$p$-tuple which is $t$ at the $k$th coordinate and $0$ elsewhere, 
\begin{align*}
\lefteqn{\bigg\| \int_{\mathbb{T}^{k-1}} \gamma_{(t_1 , 
\dots , t_{k-1}, 0, \dots , 0)} \left( \int_{\mathbb{T}}
(a - \gamma_{r_k (t_k )}(a)) K_n (t_k ) dt_k \right)}\hspace*{0.7cm} \\
\hspace*{0.5cm} &\hspace*{4cm}\times K_n (t_1 ) \cdots K_n (t_{k-1})\, 
d(t_1 , \dots , t_{k-1}) \bigg\| \\
&\leq \int_{\mathbb{T}^{k-1}}
\left\| \int_{\mathbb{T}}(a - \gamma_{r_k (t_k )}(a)) K_n (t_k ) dt_k
\right\| K_n (t_1 ) \cdots K_n (t_{k-1})\, d(t_1 , \dots , t_{k-1}) \\
&\leq \int_{\mathbb{T}} \| a - \gamma_{r_k (t_k )}(a) \| K_n (t_k )\, 
dt_k \\
&\leq L(a) \int_{\mathbb{T}} | t | K_n (t ) dt . 
\end{align*}
Estimating the integral $\int_{\mathbb{T}} |t| K_n (t)\, dt$
is a standard exercise from classical Fourier analysis
(see \cite[Exercise 3.1]{Kat}): using the fact that 
$| \sin (\pi t) | > 2 | t |$ and hence
$$ K_n (t) \leq \min\left( n+1 , \frac{1}{4(n+1)t^2} \right) $$
for $0<|t|< \frac12$,
we readily obtain, for the integral of $|t| K_n (t)$ over each of the 
intervals $[-\frac12 , -\frac{1}{2(n+1)}]$, $[-\frac{1}{2(n+1)} , 
\frac{1}{2(n+1)}]$, and $[\frac{1}{2(n+1)} , \frac12 ]$, an upper bound of 
$n^{-1}\log n$ times some constant independent of $n$, yielding the result.
\end{proof}

\begin{proposition}\label{P-torusdim}
We have
$$ \Mdim_L (A_\rho ) = p . $$
\end{proposition}

\begin{proof}
Let $\delta > 0$, and assume $\delta$ is sufficiently small so that
there is an $n\in\mathbb{N}$ such that
$$ C(n+1)^{-1}\log (n+1) \leq\delta\leq Cn^{-1}\log n $$ 
Lemma~\ref{L-convrate} then yields 
\begin{align*}
\frac{\log D(\mathcal{L}_1 , \delta )}{\log\delta^{-1}} 
&\leq \frac{\log D(\mathcal{L}_1 , 
Cn^{-1}\log n)}{\log (C(n+1)^{-1}\log (n+1))^{-1}} \\
&\leq \frac{p\log (2n+1)}{\log (Cn^{-1}\log n)^{-1}}
\end{align*}
so that
$$ \Mdim_L (A_\rho ) \leq \limsup_{n\to\infty}
\frac{p\log (2n+1)}{\log (Cn^{-1}\log n)^{-1}} = p . $$

To prove the reverse inequality, for each $n\in\mathbb{N}$ consider the set 
$$ U_n = \big\{ u_1^{k_1} u_2^{k_2} \cdots u_p^{k_p} : 
|k_i | \leq n \text{ for } i=1, \dots , p \big\} $$
of unitaries in $A_\rho$. By repeated application of the Leibniz inequality 
and using the fact that $L(u_i ) = 1$ for each $i=1, \dots , p$ we 
have the following estimate for the {\clipnorm} of an arbitrary element 
of $U_n$:
$$ L(u_1^{k_1} u_2^{k_2} \cdots u_p^{k_p}) \leq k_1 L(u_1 ) +
k_2 L(u_2 ) + \dots + k_p L(u_p ) \leq pn . $$
Thus the set $W_n = \{ (pn)^{-1}u : u\in U_n \}$ is contained in
$\mathcal{L}_1$.
Now products of distinct elements of the (self-adjoint) set $U_n$
evaluate to zero under the tracial state $\tau$, so that, in the GNS 
representation
Hilbert space associated to $\tau$ with canonical cyclic vector
$\xi_\tau$, $\pi_\tau (U_n )\xi_\tau$ forms
an orthonormal set of vectors. Thus, given $\delta > 0$  
we can apply Proposition~\ref{L-on} to obtain, for each $n\geq 1$,
$$ D(W_n , (pn)^{-1}\delta ) \geq D(U_n , \delta ) \geq
D_\tau (U_n , \delta ) \geq (1-\delta^2 )(2n+1)^p $$
so that, assuming $\delta < 1$,
\begin{align*}
\Mdim_L (A_\rho ) &\geq \limsup_{n\to\infty}\frac{\log D(W_n , 
(pn)^{-1}\delta )}{\log (pn\delta^{-1})} \\
&\geq \limsup_{n\to\infty}\frac{\log (1-\delta^2 ) + 
p\log (2n+1)}{\log (pn\delta^{-1})} \\
&= p ,
\end{align*}
as desired.
\end{proof}

\section{Product entropies}\label{S-entr}

We now study dynamics within the framework of unital $C^*$-algebras with 
Leibniz {\clipnorm}s, concentrating on iterative growth as captured in 
the ``product'' entropy of Subsection~\ref{SS-pLe} and its 
measure-theoretic version in Subsection~\ref{SS-pLes}.
That the Leibniz rule is important here can be seen by examining
the proofs of Propositions~\ref{P-Lipnumber} and \ref{P-powers} 
(although the latter only requires that $\mathcal{L}$ be closed under
multiplication).

\subsection{Product entropy}\label{SS-pLe}

We begin by introducing some notation.

\begin{notation}\label{N-entr1}
For any set $X$ we will denote by $Pf(X)$
the collection of finite subsets of $X$.
If $X_1 , X_2 , \dots , X_n$ are
subsets of the $C^*$-algebra $A$ we will use the notation 
$X_1 \cdot X_2 \cdot\,\cdots\,\cdot X_n$ or 
$\prod_{j=1}^n X_j$ to refer to the set 
$$ \{ a_1 a_2 \cdots a_n : a_i \in X_i \text{ for each } i=1, \dots ,
n \} . $$
\end{notation}

Recall from Notation~\ref{N-L} that, for a $C^*$-algebra $A$ and $r>0$, 
$A_r$ refers to the set $\{ a\in A : \| a \|\leq r \}$. For the meaning of 
$D(\cdot\, , \cdot )$ see Notation~\ref{N-dimension}.

\begin{definition}\label{D-Entp}
Let $A$ be a unital $C^*$-algebra with Leibniz {\clipnorm}
$L$, and let $\alpha\in\Aut_L (A)$. For 
$\Omega\in Pf(\mathcal{L} \cap A_1 )$ and $\delta > 0$ we define
\begin{align*}
\hspace*{5mm}\Entp_L (\alpha , \Omega , \delta ) &= 
\limsup_{n\to\infty}\frac1n \log 
D\big( \Omega\cdot\alpha(\Omega )\cdot\alpha^2 (\Omega ) \cdot\cdots\cdot
\alpha^{n-1}(\Omega ), \delta \big) , \\
\Entp_L (\alpha , \Omega ) &= \sup_{\delta > 0}
\Entp_L (\alpha , \Omega , \delta ), \\
\Entp_L (\alpha ) &= \sup_{\Omega\in Pf(\mathcal{L}\cap A_1 )}
\Entp_L (\alpha , \Omega ) .
\end{align*}
We will call $\Entp_L (\alpha )$ the {\bf product entropy} of $\alpha$.
\end{definition}

We record in the following proposition the evident fact that 
$\Entp_L (A)$ is invariant under bi-Lipschitz $^*$-isomorphisms. 

\begin{proposition} 
Let $A$ and $B$ be unital $C^*$-algebras with Leibniz
{\clipnorm}s $L_A$ and $L_B$, repectively. Let
$\alpha\in\Aut_{L_A}(A)$ and $\beta\in\Aut_{L_B}(B)$.
Suppose $\Gamma : A \to B$ is a bi-Lipschitz  
$^*$-isomorphism which intertwines $\alpha$ with $\beta$ 
(i.e., $\Gamma\circ\alpha = \beta\circ\Gamma$). Then 
$$ \Entp_L (\alpha ) = \Entp_L (\beta ) . $$
\end{proposition}

The entropy $\Entp (\alpha )$ is related to the metric dimension of $A$ by 
the following inequality, which formally parallels a familiar fact about
topological entropy (see \cite[Prop.\ 14.20]{DGS}). We remark that 
we don't know whether the Lipschitz number of a bi-Lipschitz automorphism
$\alpha$ can be strictly less than $1$, although it is evident that in
general at least one of $\alpha$ and $\alpha^{-1}$ must have Lipschitz
number at least $1$. 

\begin{proposition}\label{P-Lipnumber}
If $\alpha\in\Aut_L (A)$ and $\Mdim_L (A)$ is finite then
$$ \Entp_L (\alpha ) \leq \Mdim_L (A)\cdot\log\max(\lambda , 1) $$
where $\lambda$ is the Lipschitz number of $\alpha$. 
\end{proposition}

\begin{proof}
Let $\Omega\in Pf(\mathcal{L}\cap A_1 , \delta )$ and $\delta > 0$.
Set $M = \max_{a\in\Omega} L(a)$. Then by repeated application of the
Leibniz inequality we see that elements of the set 
$$ \Omega_n = \Omega\cdot\alpha(\Omega )\cdot\alpha^2 (\Omega ) 
\cdot\cdots\cdot\alpha^{n-1}(\Omega ) $$
have {\clipnorm} at most $M(1 + \lambda + \lambda^2 + \cdots + \lambda^{n-1})$,
which is bounded above by $Mn\lambda^n$. Hence $\mathcal{L}_1$
contains the set $\{ (Mn\lambda^n )^{-1}a : a\in\Omega_n \}$, which we will 
denote simply by $(Mn\lambda^n )^{-1}\Omega_n$. It follows that
\begin{align*}
\Entp_L (\alpha , \Omega , \delta ) &= 
\limsup_{n\to\infty}\frac1n \log 
D( \Omega_n , \delta ) \\
&= \limsup_{n\to\infty}\frac1n \log D( (Mn\lambda^n )^{-1}\Omega_n , 
(Mn\lambda^n )^{-1}\delta ) \\
&\leq \limsup_{n\to\infty}\frac1n
\log D(\mathcal{L}_1 , (Mn\lambda^n )^{-1}\delta ) . 
\end{align*}
If $\lambda < 1$ then this last limit supremum is clearly zero. 
If on the other hand $\lambda \geq 1$ then
\begin{align*}
\lefteqn{\limsup_{n\to\infty}\frac1n
\log D(\mathcal{L}_1 , (Mn\lambda^n )^{-1}\delta )}\hspace*{2.4cm} \\ 
\hspace*{2.3cm} &\leq \limsup_{n\to\infty}\frac1n
\frac{\log D(\mathcal{L}_1 , (Mn\lambda^n )^{-1}\delta )}{\log
(Mn\lambda^n \delta^{-1})}\log (Mn\lambda^n \delta^{-1}) \\
&= \limsup_{n\to\infty}\frac{\log D(\mathcal{L}_1 , 
(Mn\lambda^n )^{-1}\delta )}{\log (Mn\lambda^n \delta^{-1})}\cdot
\lim_{n\to\infty}\frac1n \log (Mn\lambda^n \delta^{-1}) \\
&= \Mdim_L (A)\cdot\log\lambda .
\end{align*}
We thus obtain the result by taking the supremum over all $\Omega$
and $\delta$.
\end{proof}

\begin{corollary}\label{C-isomentropy}
If $\Mdim_L (A)$ is finite and $\alpha\in\Aut_L (A)$ is Lipschitz 
isometric then $\Entp_L (\alpha ) = 0$. In particular $\Entp_L (\id_A ) = 0$.
\end{corollary}

Corollary~\ref{C-isomentropy} shows that the appropriate domain for
our notion of entropy as a measure of dynamical growth is the
class of {\clipnorm}ed unital $C^*$-algebras $A$ for which $\Mdim_L (A)$ is
finite, in analogy to the situation of topological approximation entropies 
\cite{Br,Voi} which function under conditions of ``finiteness'' 
like nuclearity or exactness.

\begin{proposition}\label{P-powers}
If $\alpha\in\Aut_L (A)$ and $k\in\mathbb{Z}$ then 
$\Entp_L (\alpha^k ) = |k|\, \Entp_L (\alpha )$. 
\end{proposition}

\begin{proof}
Suppose first that $k\geq 0$. Let $\Omega\in Pf(\mathcal{L} \cap A_1 )$ and 
$\delta > 0$, and suppose $1\in\Omega$. Then 
$$ \prod_{j=0}^{n-1}\alpha^{jk}(\Omega ) \subset 
\prod_{j=0}^{(n-1)k}\alpha^{j}(\Omega ) $$
so that
\begin{align*}
\Entp_L (\alpha^k , \Omega , \delta ) &=
\limsup_{n\to\infty}\frac1n \log 
D\Bigg( \prod_{j=0}^{n-1}\alpha^{jk}(\Omega ) , \delta \Bigg) \\
&\leq k \limsup_{n\to\infty}\frac{1}{kn} \log 
D\Bigg( \prod_{j=0}^{(n-1)k}\alpha^{j}(\Omega ) , \delta \Bigg) \\
&= k\, \Entp_L (\alpha , \Omega , \delta ) .
\end{align*}
On the other hand setting $\Omega_k = \prod_{j=0}^{k-1}\alpha^{j}(\Omega )$
(which is contained in $Pf(\mathcal{L}\cap A_1 )$ in view of the 
Leibniz rule) we have
$$ \prod_{j=0}^{\lfloor \frac{n}{k} \rfloor}\alpha^{j}(\Omega_k ) 
\subset\prod_{j=0}^{n-1}\alpha^j (\Omega ) $$
so that
\begin{align*}
\Entp_L (\alpha^k , \Omega_k , \delta ) &=
\limsup_{n\to\infty}\frac{k}{n} \log 
D\Bigg( \prod_{j=0}^{\lfloor \frac{n}{k} \rfloor}\alpha^j (\Omega_k ) , 
\delta \Bigg) \\
&\leq k \limsup_{n\to\infty}\frac1n \log 
D\Bigg( \prod_{j=0}^{n-1}\alpha^{j}(\Omega ) , 
\delta \Bigg) \\
&= k\, \Entp_L (\alpha , \Omega , \delta ) .
\end{align*}
and hence
$$ \Entp_L (\alpha^k , \Omega_k , \delta ) \leq
k\, \Entp_L (\alpha , \Omega , \delta ) . $$
Taking the supremum over all $\Omega\in Pf(\mathcal{L} \cap A_1 )$ 
and $\delta > 0$ yields $\Entp_L (\alpha^k ) = k\, \Entp_L (\alpha )$.

To prove the assertion for $k<0$ it suffices, in view of the first part, 
to show that $\Entp_L (\alpha^{-1}) = \Entp_L (\alpha )$. Since  
$$ \alpha^{-n+1}\Bigg(\prod_{j=0}^{n-1}\alpha^j (\Omega )\Bigg) =
\prod_{j=0}^{n-1}\alpha^{-j} (\Omega ) $$
we have
$$ D\Bigg( \prod_{j=0}^{n-1}\alpha^j (\Omega ) , \delta \Bigg) = 
D\Bigg( \prod_{j=0}^{n-1}\alpha^{-j}(\Omega ) , \delta \Bigg) $$
and hence
$$ \Entp_L (\alpha , \Omega , \delta ) =
\Entp_L (\alpha^{-1}, \Omega , \delta ) , $$
from which we reach the conclusion by taking the supremum over all 
$\Omega\in Pf(\mathcal{L} \cap A_1 )$ and $\delta > 0$.
\end{proof}

The following proposition is clear from Definition~\ref{D-Entp}.

\begin{proposition}
Let $A$ be a unital $C^*$-algebra with {\clipnorm} $L_A$ and
$B\subset A$ a unital $C^*$-subalgebra with {\clipnorm} $L_B$ such that
$L_B$ is the restriction of $L_A$ to $B$. Suppose that there is a
$C^*$-norm contractive idempotent linear map of $A$ onto $B$.
If $\alpha\in\Aut_L (A)$ leaves $B$ invariant then
$$ \Entp_{L_B}(\alpha |_B ) \leq \Entp_{L_A}(\alpha ) . $$
\end{proposition}

\begin{proposition}
Let $A$ and $B$ be unital $C^*$-algebras, $L_A$ a Leibniz
{\clipnorm} on $A$, $\phi : A \to B$ a surjective unital 
$^*$-homomorphism, and $L_B$ the Leibniz {\clipnorm} induced on 
$B$ via $\phi$. Suppose there exists a positive
$C^*$-norm contractive (not necessarily unital) Lipschitz map
$\psi : B \to A$ such that $\phi\circ\psi = \id_B$. Let 
$\alpha\in\Aut_{L_A}(A)$ and $\beta\in\Aut_{L_B}(B)$ 
and suppose $\phi\circ\alpha = \beta\circ\phi$. Then
$$ \Entp_{L_B}(\beta ) \leq \Entp_{L_A}(\alpha ) . $$
\end{proposition}

\begin{proof}
Let $\Omega\in Pf(\mathcal{L}^B \cap B_1 )$ and $\delta > 0$. Since
$\psi$ is norm-decreasing we have $\psi (\Omega ) \in 
Pf(\mathcal{L}^A \cap A_1 )$. Now if $X\in\mathcal{F}(A)$ is such that
$$ \psi (\Omega ) \cdot\alpha(\psi (\Omega ))\cdot\cdots\cdot
\alpha^{n-1}(\psi (\Omega )) \subset_\delta X $$
then
\begin{align*}
\Omega\cdot\beta (\Omega ) \cdot\cdots\cdot\beta^{n-1}(\Omega )
&= (\phi\circ\psi )(\Omega ) \cdot\beta ((\phi\circ\psi )(\Omega ))
\cdot\cdots\cdot\beta^{n-1}((\phi\circ\psi )(\Omega )) \\
&= \phi (\psi (\Omega ))\cdot\phi (\alpha (\psi (\Omega )))\cdot\cdots\cdot
\phi (\alpha^{n-1}(\psi (\Omega ))) \\
&= \phi (\psi (\Omega ) \cdot\alpha(\psi (\Omega ))\cdot\cdots\cdot
\alpha^{n-1}(\psi (\Omega ))) \\
&\subset_\delta \phi (X)
\end{align*}
and so 
$$ D(\Omega\cdot\beta (\Omega ) \cdot\cdots\cdot\beta^{n-1}(\Omega ), 
\delta )\leq D(\psi (\Omega ) \cdot\alpha(\psi (\Omega ))\cdot\cdots\cdot
\alpha^{n-1}(\psi (\Omega )), \delta ) , $$
from which the proposition follows.
\end{proof}

\subsection{Product entropy with respect to an invariant 
state}\label{SS-pLes}

We define now a version of $\Mdim_L (A)$ relative to a dynamically
invariant state $\sigma$. As in Subsection~\ref{SS-pLe} we are
assuming that $L$ is a Leibniz {\clipnorm}. For the meaning of
$D_\sigma (\cdot\, , \cdot )$ see Notation~\ref{N-dimension}.

\begin{definition}
Let $\alpha\in\Aut_L (A)$ and let $\sigma$ be a state of $A$
which is $\alpha$-invariant, i.e., $\sigma\circ\alpha = \sigma$. For 
$\Omega\in Pf(\mathcal{L} \cap A_1 )$ and $\delta > 0$ we define
\begin{align*}
\hspace*{5mm}\Entp_{L, \sigma}(\alpha , \Omega , \delta ) &= 
\limsup_{n\to\infty}\frac1n \log 
D_\sigma \big(\Omega\cdot\alpha(\Omega )\cdot\alpha^2 (\Omega ) 
\cdot\cdots\cdot\alpha^{n-1}(\Omega ), \delta \big) , \\
\Entp_{L,\sigma} (\alpha , \Omega ) &= \sup_{\delta > 0}
\Entp_{L,\sigma} (\alpha , \Omega , \delta ), \\
\Entp_{L,\sigma} (\alpha ) &= \sup_{\Omega\in Pf
(\mathcal{L} \cap A_1 )} \Entp_{L, \sigma}(\alpha , \Omega ) .
\end{align*}
We will call $\Entp_{L,\sigma} (\alpha )$ the 
{\bf product entropy} of $\alpha$ with respect to 
$\sigma$.
\end{definition}

The following two propositions follow immediately from the definition. 

\begin{proposition} 
Let $A$ and $B$ be unital $C^*$-algebras with respective Leibniz
{\clipnorm}s $L_A$ and $L_B$. Let
$\alpha\in\Aut_{L^A}(A)$ and $\beta\in\Aut_{L^B}(B)$, 
and let $\sigma$ and $\omega$ 
be $\alpha$- and $\beta$-invariant states on $A$ and $B$, respectively. 
Suppose $\Gamma : A \to B$ is a bi-Lipschitz 
$^*$-isomorphism such that $\Gamma\circ\alpha = \beta\circ\Gamma$
and $\omega\circ\Gamma = \sigma$. Then 
$$ \Entp_{L, \sigma}(\alpha ) = \Entp_{L, \omega}(\beta ) . $$
\end{proposition}

\begin{proposition}
Let $A$ be a unital $C^*$-algebra with Leibniz {\clipnorm} $L_A$ and
$B\subset A$ a unital $C^*$-subalgebra with Leibniz {\clipnorm} $L_B$ such 
that $L_B$ is the restriction of $L_A$ to $B$. Let $\sigma$ be a state on
$A$ with $\sigma\circ\alpha = \sigma$, and suppose that there is a
idempotent linear map of $A$ onto $B$ which is contractive for the
Hilbert space norm under the GNS construction associated to $\sigma$. If 
$\alpha\in\Aut_L (A)$ leaves $B$ invariant then
$$ \Entp_{L_B , \sigma}(\alpha |_B ) \leq \Entp_{L_A , \sigma}(\alpha ) . $$
\end{proposition}

The next proposition can be established in the same way as its
counterpart Proposition~\ref{P-powers} in Subsection~\ref{SS-pLe}.

\begin{proposition}
If $\alpha\in\Aut_L (A)$, $\sigma$ is an $\alpha$-invariant state
on $A$, and $k\in\mathbb{Z}$, then 
$\Entp_{L,\sigma} (\alpha^k ) = |k|\, \Entp_{L,\sigma} (\alpha )$. 
\end{proposition}

\begin{proposition}\label{P-comparison}
If $\alpha\in\Aut_L (A)$ and $\sigma$ is an $\alpha$-invariant 
state on $A$ then
$$ \Entp_{L,\sigma} (\alpha ) \leq \Entp_L (\alpha ) . $$
\end{proposition}

\begin{proof}
It suffices to show that, for a given $\Omega\in Pf(\mathcal{L}\cap A_1 )$ 
and $\delta > 0$, 
$$ D_\sigma (\Omega , \delta ) \leq D(\Omega  , \delta ) , $$
and for this inequality we need only observe that if $X$ is a 
finite-dimensional subspace of $A$ such that $\Omega\subset_\delta X$,
then whenever $a\in\Omega$ and $x\in X$ satisfy $\| a - x \| < \delta$
we have
$$ \| \pi (a)\xi_\sigma - \pi (x)\xi_\sigma \|_\sigma =
\| \pi (a - x)\xi_\sigma \|_\sigma \leq \| \pi (a - x) \| \leq \| a - x \|
< \delta , $$
so that $\pi (X)\xi_\sigma$ is a subspace of $\mathcal{H}_\sigma$ with
$\pi (\Omega )\xi_\sigma \subset_\delta \pi (X)\xi_\sigma$ and
$\dim \pi (\Omega )\xi_\sigma \leq \dim X$.
\end{proof}

\begin{corollary}\label{C-isomstateentropy}
If $\Mdim_L (A)$ is finite and $\alpha\in\Aut_L (A)$ is Lipschitz 
isometric then $\Entp_{L,\sigma}(\alpha ) = 0$. In particular 
$\Entp_{L,\sigma} (\id_A ) = 0$.
\end{corollary}

\begin{proof}
This follows by combining Proposition~\ref{P-comparison} with 
Corollary~\ref{C-isomentropy}.
\end{proof}

\section{Tensor product shifts}\label{S-shift} 

The fundamental prototypical system for topological entropy is the
shift on the infinite product
$\{ 1, \dots , p \}^{\mathbb{Z}}$, with entropy $\log p$.
Here we consider the noncommutative analogue of this map, 
the (right) shift on the infinite tensor product $M_p^{\otimes\mathbb{Z}}$
of $p\times p$ matrix algebras $M_p$ over $\mathbb{C}$, here with
the Leibniz {\clipnorm} $L$ furnished by the infinite product 
$\gamma^{\otimes\mathbb{Z}} : G^{\mathbb{Z}} \to
\Aut \big( M_p^{\otimes\mathbb{Z}}\big)$ of Weyl actions
and length function $\ell_\lambda$ (for a given $0< \lambda < 1$)
as described in Subsection~\ref{SS-UHF}.

Before computing the entropy of the shift we will show that it is a
bi-Lipschitz $^*$-automorphism.

\begin{proposition}\label{P-shiftLip}
The shift $\alpha$ on $M_p^{\otimes\mathbb{Z}}$ is a bi-Lipschitz
$^*$-automorphism, and $\alpha$ and its inverse have 
Lipschitz numbers bounded by $\lambda$.
\end{proposition}

\begin{proof}
Let $T : G^{\mathbb{Z}} \to G^{\mathbb{Z}}$ be the 
right shift homeomorphism. Then it is readily seen that if $a$ is an 
elementary tensor in $M_p^{[m,n]}\subset M_p^{\otimes\mathbb{Z}}$ for some
$m,n\in\mathbb{Z}$ then $\gamma_g^{\otimes\mathbb{Z}} 
(\alpha (a)) = \alpha \big( \gamma_{Tg}^{\otimes\mathbb{Z}}(a)\big)$
for all $g\in G^{\mathbb{Z}}$,
and since such $a$ generate $M_p^{\otimes\mathbb{Z}}$ we have
$\gamma_g^{\otimes\mathbb{Z}} \circ\alpha = 
\alpha\circ\gamma_{Tg}^{\otimes\mathbb{Z}}$ for all $g\in
G^{\mathbb{Z}}$. Thus, for any $a\in M_p^{\otimes\mathbb{Z}}$,
\begin{align*}
L(\alpha (a)) &= \sup_{g\in G^{\mathbb{Z}}\setminus\{ e \}}
\frac{\big\| \gamma_g^{\otimes\mathbb{Z}} (\alpha (a)) - 
\alpha (a) \big\|}{\ell_\lambda (g)} \\
&= \sup_{g\in G^{\mathbb{Z}}\setminus\{ e \}}
\frac{\big\| \alpha \big( \gamma_{Tg}^{\otimes\mathbb{Z}}(a)\big) - 
\alpha (a) \big\|}{\ell_\lambda (g)} \\
&\leq \sup_{g\in G^{\mathbb{Z}}\setminus\{ e \}}
\frac{\big\| \gamma_{Tg}^{\otimes\mathbb{Z}}(a) - a \big\|}{\ell_\lambda (Tg)} 
\cdot \sup_{g\in G^{\mathbb{Z}}\setminus\{ e \}}
\frac{\ell_\lambda (Tg)}{\ell_\lambda (g)} \\
&\leq L(a) \cdot L(T) ,
\end{align*}
where $L(T)$ is the Lipschitz number of the homeomorphism $T$ with respect
to the metric defining $\ell_\lambda$ (see Subsection~\ref{SS-UHF}), and it 
is straightforward to verify that $L(T) = \lambda$. We can argue similarly
for $\alpha^{-1}$ to reach the desired conclusion. 
\end{proof} 

\begin{proposition}\label{P-shiftlowerbd}
Let $\alpha$ be the shift on $M_p^{\otimes\mathbb{Z}}$ and 
$\tau = \text{tr}_p^{\otimes\mathbb{Z}}$ the
unique (and hence $\alpha$-invariant) tracial state on 
$M_p^{\otimes\mathbb{Z}}$. Then 
$$ \Entp_{L , \tau}(\alpha ) \geq 2\log p . $$
\end{proposition}

\begin{proof}
Let $u,v\in M_p^{\otimes\mathbb{Z}}$ be the Weyl generators for the 
zeroeth copy of $M_p$ (identified as a subalgebra of 
$M_p^{\otimes\mathbb{Z}}$) and let $\Omega$ be the finite subset
$\{ u^i v^j : 0\leq i,j \leq k-1 \}$ of 
$\mathcal{L}\cap \big( M_p^{\otimes\mathbb{Z}}\big)_1$. 
Then the set $\Omega_n = \Omega\cdot\alpha(\Omega )\cdot\alpha^2 (\Omega ) 
\cdot\cdots\cdot\alpha^{n-1}(\Omega )$ is precisely the subset
$$ \{ u^{i_0}v^{j_0} \otimes u^{i_1}v^{j_1} \otimes\cdots\otimes
u^{i_{n-1}}v^{j_{n-1}} : 0\leq i_k ,j_k \leq p-1 \text{ for }
k=0, \dots , n-1 \} $$
of $M_p^{\otimes [0,n]}$ as considered sitting in 
$M_p^{\otimes\mathbb{Z}}$. Thus $\pi_\tau (\Omega_n )\xi_\tau$
is an orthonormal set in the GNS representation Hilbert space 
associated to $\tau$ with canonical cyclic vector $\xi_\tau$
(see the second half of the proof of Proposition~\ref{P-Mpdim}), and
so by Lemma~\ref{L-on} for any $\delta > 0$ we have
$$ D_\tau (\Omega_n , \delta ) \geq (1-\delta^2 ) \cardn 
(\pi_\tau (\Omega_n )\xi_\tau ) = (1-\delta^2 ) p^{2n} . $$
Thus if $\delta < 1$ we obtain
$$ \Entp_{L, \sigma}(\alpha , \Omega , \delta ) = 
\limsup_{n\to\infty}\frac1n \log D_\sigma (\Omega_n , \delta ) \geq
2\log p , $$
which yields the proposition.
\end{proof}

Note that by Propositions~\ref{P-Lipnumber}, \ref{P-Mpdim}, and
\ref{P-shiftLip} the shift $\alpha$ satisfies 
$$ \Entp_L (\alpha ) \leq 4\log p . $$
The following proposition yields the sharp upper bound of $2\log p$.
 
\begin{proposition}\label{P-shiftupperbd}
With $\alpha$ the shift we have 
$$ \Entp_L (\alpha ) \leq 2\log p . $$
\end{proposition}

\begin{proof}
Let $\Omega\in Pf\big( \mathcal{L}\cap 
\big( M_p^{\otimes\mathbb{Z}}\big)_1 \big)$ and $\delta > 0$.
Set $C = \max_{a\in\Omega} L_\lambda (a)$. For each $n$ consider the 
conditional expectation 
$E_n : M_p^{\otimes\mathbb{Z}} \to M_p^{\otimes [-n,n]}$ given by
$$ E_n (a) = \int_{G^{\mathbb{Z}\setminus [-n,n]}} 
\gamma_{g}^{\otimes\mathbb{Z}}(a)\, dg , $$
where $dg$ is normalized Haar measure on $G^{\mathbb{Z}}$. We then have
\begin{align*}
\| E_n (a) - a \| &= \bigg\| \int_{G^{\mathbb{Z}\setminus [-n,n]}} 
\big( \gamma_{g}^{\otimes\mathbb{Z}}(a) - a\big) \, dg \bigg\| \\ 
&\leq \int_{G^{\mathbb{Z}\setminus [-n,n]}}
\big\| \gamma_{g}^{\otimes\mathbb{Z}}(a) - a \big\| \, dg \\
&\leq \int_{G^{\mathbb{Z}\setminus [-n,n]}} C\ell_\lambda (g)
\, dg \\
&\leq \frac{2C\lambda^{n+1}}{1-\lambda} .
\end{align*}
If $a_1 \dots , a_n \in\Omega$ then, estimating the norm of
differences of products in the usual way and using the fact that 
the conditional expectations are norm-decreasing, we have
\begin{align*}
\lefteqn{\big\| E_{\lceil\sqrt{n}\rceil}(a_1 ) \alpha 
(E_{\lceil\sqrt{n}\rceil}(a_2 )) \cdots\alpha^{n-1} 
(E_{\lceil\sqrt{n}\rceil}(a_n )) - a_1 \alpha (a_2 )
\cdots\alpha^{n-1}(a_n ) \big\|}\hspace*{6cm} \\
\hspace*{5.7cm} &\leq \sum_{k=1}^n \big\| 
\alpha^{k-1}(E_{\lceil\sqrt{n}\rceil}(a_k )) -
\alpha^{k-1}(a_k ) \big\| \\
&= \sum_{k=1}^n \big\| E_{\lceil\sqrt{n}\rceil}(a_k ) - a_k \big\| \\
&\leq \frac{2Cn\lambda^{\lceil\sqrt{n}\rceil + 1}}{1-\lambda} , 
\end{align*}
which is smaller than $\delta$ for all $n$ greater than some
$n_0 \in \mathbb{N}$ (here $\lceil\cdot\rceil$ denotes the ceiling
function).

Next we observe that the product
$$ E_{\lceil\sqrt{n}\rceil}(a_1 ) \alpha (E_{\lceil\sqrt{n}\rceil}
(a_2 )) \cdots\alpha^{n-1}(E_{\lceil\sqrt{n}\rceil}(a_n )) $$ 
is contained in the
subalgebra $M_p^{\otimes[-\lceil\sqrt{n} \rceil , \lceil\sqrt{n} 
\rceil + n]}$ of $M_p^{\otimes\mathbb{Z}}$, and this subalgebra
has linear dimension $p^{2(2\lceil\sqrt{n}\rceil + n)}$. In view of the
first paragraph, for all $n\geq n_0$ the set $\Omega_n$ is
approximately contained in $M_p^{\otimes[-\lceil\sqrt{n} \rceil , 
\lceil\sqrt{n} \rceil + n]}$ to within $\delta$, and so we
have
$$ \Entp_L (\alpha , \Omega , 2\delta )
\leq\limsup_{n\to\infty} \frac1n \log p^{2(2\lceil\sqrt{n}\rceil + n)} = 
2\log p . $$
The proposition now follows by taking the supremum over all $\Omega$
and $\delta$.
\end{proof} 

As a consequence of Propositions \ref{P-shiftlowerbd},
\ref{P-shiftupperbd}, and \ref{P-comparison} we obtain the following.

\begin{proposition}
With $\alpha$ the shift and $\tau$ the unique tracial state on 
$M_p^{\otimes\mathbb{Z}}$ we have
$$ \Entp_L (\alpha ) =
\Entp_{L , \tau}(\alpha ) =  2\log p . $$
\end{proposition}

\section{Noncommutative toral automorphisms}\label{S-tori}

Let $A_\rho$ be a noncommutative $p$-torus with generators
$u_1 , \dots , u_p$, canonical ergodic action $\gamma : \mathbb{T}^p
\to\Aut (A_\rho )$, and associated Leibniz {\clipnorm} $L$ and  
$\gamma$-invariant tracial state $\tau$, 
as defined in Subsection~\ref{SS-noncommtorusdim}.
We let $\pi_\tau : A_\rho \to\mathcal{B}(\mathcal{H}_\tau )$ be the
GNS representation associated to $\tau$, with canonical
cyclic vector $\xi_\tau$.
Let $T = (s_{ij})$ be a $p\times p$ integral matrix with 
$\det T = \pm 1$, and suppose that $T$ defines an automorphism 
$\alpha_T$ of $A_\rho$ via the specifications
$$ \alpha_T (u_j ) = u_1^{s_{1j}} \cdots u_p^{s_{pj}} $$  
on the generators (this will always be the case if $\det T = 1$ owing to 
the universal property of noncommutative tori).
These noncommutative versions of toral automorphisms
were introduced in the case $p=2$ in \cite{Wat} and \cite{Bre}.
Since $\tau$ is zero on products of powers of generators which are not
equal to the unit, we see that it is invariant under the automorphism 
$\alpha_T$ and the action $\gamma$. 
Fix a $t = (t_1 , \dots , t_p )\in\mathbb{T}^p \cong (\mathbb{R} / 
\mathbb{Z})^p$ and consider the automorphism $\gamma_t$ coming from
the action $\gamma$.
We will compute the entropies $\Entp_L (\alpha_T \circ\gamma_t )$ 
and $\Entp_{L , \tau}(\alpha_T \circ\gamma_t )$ and 
furthermore show that their common value bounds above the entropies 
$\Entp_L (\Ad u \circ\alpha_T \circ\gamma_t )$ and
$\Entp_{L, \tau}(\Ad u \circ\alpha_T \circ\gamma_t )$ for any unitary 
$u\in\mathcal{L}$. We remark that in the case $p=2$, 
when $A_\rho$ is a rotation $C^*$-algebra $A_\theta$, Elliott showed 
in \cite{Ell86} that if the angle $\theta$ satisfies a generic Diophantine 
property then all automorphisms preserving the dense
$^*$-subalgebra of smooth elements (i.e., all ``diffeomorphisms'')
are of the form $\Ad u \circ\alpha_T \circ\gamma_t$ where $u$ is
a smooth unitary (and hence of finite {\clipnorm}).   

\begin{proposition}\label{P-torusLip}
The $^*$-automorphism $\alpha = \Ad u \circ\alpha_T \circ\gamma_t$ is 
bi-Lipschitz, and $\alpha$ and its 
inverse have Lipschitz numbers bounded by 
$$ 2\| T \|(1+2L(u)\diam (S(A))) $$ 
and 
$$ 2\| T^{-1} \|(1+2L(u)\diam (S(A))) , $$ 
respectively, where $\| T \|$ and $\| T^{-1} \|$ are the respective norms
of $T$ and $T^{-1}$ as operators on the real inner product space 
$\mathbb{R}^p$.
\end{proposition}

\begin{proof}
If we consider $T$ as acting on $\mathbb{T}^p$ then 
$\gamma_g \circ\alpha = \alpha\circ\gamma_{Tg}$ for all $g\in\mathbb{T}^p$,
as can be seen by checking this equation on the generators 
$u_1 , \dots , u_p$. As in the proof of Proposition~\ref{P-shiftLip}
we thus have, for any $a\in\mathcal{L}$, the bound
$$ L(\alpha (a)) \leq L(a) \cdot L(T) $$
where $L(T)$ is the Lipschitz number of the homeomorphism $T$. If we
consider $T$ as an operator on $\mathbb{R}^p$,
then its Lipschitz number is $\| T \|$ by definition of the operator norm, 
and so by linearity the Lipschitz number $L(T)$ of $T$ on the quotient 
$\mathbb{T}^p \cong \mathbb{R}^p / \mathbb{Z}^p$ must 
again be $\| T \|$. Next note that 
$\gamma_t$ is isometric, for if $a\in\mathcal{L}$ then
$$ L(\gamma_t (a)) = \sup_{s\in\mathbb{T}^p \setminus\{ 0 \}}
\frac{\| \gamma_{s+t}(a) - \gamma_t (a) \|}{\ell 
(s)} = \sup_{s\in\mathbb{T}^p \setminus\{ 0 \}}
\frac{\| \gamma_s (a) - a \|}{\ell (s)} = L(a) . $$
Also, since $L$ is readily checked to be lower semicontinuous,
by Proposition~\ref{P-inner} the Lipschitz number of
$\Ad u$ is bounded by $2(1+2L(u)\diam (S(A)))$. Thus by 
Proposition~\ref{P-composition} we get the desired bound on the Lipschitz
number of $\Ad u \circ\alpha_T \circ\gamma_t$. 
A similar argument can be applied to 
$(\Ad u \circ\alpha_T \circ\gamma_t )^{-1} = \gamma_{-t} \circ
\alpha_{T^{-1}}\circ\Ad u^*$.  
\end{proof}

\begin{proposition}\label{P-lowerbound}
We have
$$ \Entp_{L, \tau} (\alpha_T \circ\gamma_t ) \geq
\sum_{|\lambda_i | \geq 1} \log |\lambda_i | $$
where $\lambda_1 , \cdots , \lambda_p$
are the eigenvalues of $T$ counted with spectral multiplicity.
\end{proposition}

\begin{proof}
Let $K$ be a finite subset of $\mathbb{Z}^p$ and set
$$ U_K = \big\{ u_1^{k_1} u_2^{k_2} \cdots u_p^{k_p} : 
(k_1 , \dots , k_p ) \in K \big\} . $$
The elements of $U_K$, being products of powers of generators, all have
finite {\clipnorm}. 
Observe that $\alpha_T \circ\gamma_t$ takes a product of the form 
$\eta u_1^{k_1} u_2^{k_2}\cdots u_p^{k_p}$, with $\eta$ a complex
number of modulus one, to a product of the same form, with the
exponents on the $u_i$'s respecting the action of the group automorphism
$\zeta_T$ of $\mathbb{Z}^p$ defined via the action of $T$. 
Thus if $K$ is a finite subset of $\mathbb{Z}^p$ then the set
$$ U_K \cdot (\alpha_T \circ\gamma_t )(U_K )\cdot\cdots\cdot 
(\alpha_T \circ\gamma_t )^{n-1}(U_K )$$
contains a subset $U_{K,n}$ of the form
$$ \big\{ \eta_{(k_1 , \dots , k_p )}u_1^{k_1} u_2^{k_2}\cdots u_p^{k_p}
: (k_1 , \dots , k_p )\in K + \zeta_T K + \cdots + \zeta_T^{n-1}K \big\} $$
where each $\eta_{(k_1 , \dots , k_p )}$ is a complex number of modulus
one. Note that $\pi_\tau (U_{K,n})\xi_\tau$ is an orthonormal
set of vectors in the GNS representation Hilbert space 
associated to $\tau$ with canonical cyclic vector $\xi_\tau$, 
since the product of any two distinct
vectors in this set is a scalar multiple of a product of the form
$u_1^{k_1} u_2^{k_2}\cdots u_p^{k_p}$ with the $k_i$'s not all zero, in which 
case evaluation under $\tau$ yields zero. It thus follows from 
Lemma~\ref{L-on} that if $\delta > 0$ then
\begin{align*}
D_\tau (U_{K,n}, \delta ) &\geq (1-\delta^2 )
\cardn (\pi (U_{K,n})\xi_\tau ) \\
&= (1-\delta^2 ) \cardn (K + \zeta_T K + \cdots + \zeta_T^{n-1}K ) , 
\end{align*}
so that whenever $\delta < 1$ we get
\begin{align*}
\Entp_{L,\sigma} (\alpha_t \circ\gamma_t , U_K , \delta ) &= 
\limsup_{n\to\infty}\frac1n \log 
D_\tau (U_K \cdot\alpha (U_K )\cdot\cdots\cdot\alpha^{n-1}
(U_K ), \delta ), \\
&\geq \limsup_{n\to\infty}\frac1n \log D_\tau (U_{K,n}, \delta ) \\
&\geq \limsup_{n\to\infty} \frac1n \log \cardn (K + \zeta_T K + 
\cdots + \zeta_T^{n-1}K ) .
\end{align*}
We thus reach the desired conclusion by recalling from the computation 
of the discrete Abelian group entropy of $\zeta_T$ \cite{Pet} that 
$$ \lim_K \limsup_{n\to\infty} \frac1n \log \cardn (K + \zeta_T K + \cdots
+ \zeta_T^{n-1}K ) = \sum_{|\lambda_i | \geq 1} 
\log |\lambda_i | $$
where the limit is taken with respect to the net of finite subsets $K$ of 
$\mathbb{Z}^p$.
\end{proof}

To compute upper bounds we need a couple of lemmas.

\begin{lemma}\label{L-card}
Let $\zeta_T$ be the group automorphism of $\mathbb{Z}^p$ defined via
the action of an $p\times p$ integral matrix $T$ with $\det (T) =\pm 1$. 
Let $\lambda_1 , \cdots , \lambda_p$
be the eigenvalues of $T$ counted with spectral multiplicity. For each
$m\in\mathbb{N}$ let $K_m$ be the cube
$$ \{ (k_1 , \dots , k_p ) \in\mathbb{Z}^p : |k_i | \leq m \text{ for each } 
i=1, \dots , p \}  $$
and define recursively for $n\geq 0$ the 
sets $L_{m,n}\in\mathbb{Z}^p$
by $L_{m,0} = K_m$ and 
$$ L_{m, n+1} = \zeta_T (L_{m,n}) + K_m . $$ 
Then for every $\delta > 0$ there is a $Q>0$ such that, for all 
$m, n\in\mathbb{N}$,
$$ \cardn (L_{m,0} + L_{m,1} + \cdots + L_{m,n-1}) 
\leq Q(mn^2)^p (1+\delta )^n \prod_{|\lambda_i | \geq 1}|\lambda_i |^n . $$
\end{lemma}  

\begin{proof}
For any subset $K$ of $\mathbb{Z}^p$ we will denote its convex hull
as a subset of $\mathbb{R}^p$ by $\tilde{K}$. With $\zeta_T$ 
also referring to the linear map on $\mathbb{R}^p$ defined by $T$, 
we consider the convex set $\tilde{L}_{m,0} + \tilde{L}_{m,1} + \cdots 
+ \tilde{L}_{m,n-1}$. By amplifying this
set by a linear factor of $2^p$ we can ensure that it contains
every cube of unit side length centred at some point in $L_{m,0} + 
L_{m,1} + \cdots + L_{m,n-1}$, so that
$$ \cardn(L_{m,0} + L_{m,1} + \cdots + L_{m,n-1})
\leq 2^p \vol (\tilde{L}_{m,0} + \tilde{L}_{m,1} + \cdots 
+ \tilde{L}_{m,n-1}) . $$
To estimate this volume on the right we assemble a basis $\mathcal{B}$
of $\mathbb{R}^p$ by picking a basis for the spectral subspace
associated to each real eigenvalue and each pair of conjugate complex
eigenvalues. Working from this point on with respect to the basis
$\mathcal{B}$, we note that the sets $\tilde{K}_m$ are now
parallelipipeds, and they can be contained in cubes $B_m$ centred at
$0$ of side length $rm$ for some $r>0$ independent of $m$ by the linearity 
of our basis change. If we define the sets $M_{m,n}$ recursively
by $M_{m,0} = B_m$ and
$$ M_{m,n+1} = \zeta_T (M_{m,n}) + B_m $$
then the set $M_{m,0} + M_{m,1} + \cdots + M_{m,n-1}$ is a $p$-dimensional
rectangular box which is centred at the origin with each face
perpendicular to some coordinate axis, and this box contains
$\tilde{L}_{m,0} + \tilde{L}_{m,1} + \cdots + \tilde{L}_{m,n-1}$,
so that it suffices to show that 
$$ \vol (M_{m,0} + M_{m,1} + \cdots + M_{m,n-1}) $$
is bounded by the last expression in the lemma statement for some $C>0$.

Let $v$ be a vector in $\mathcal{B}$ associated to a real
eigenvalue $\lambda$ or a complex conjugate pair $\{ \lambda , 
\bar{\lambda} \}$. We can then find a $Q>0$ such that for all $n\in
\mathbb{N}$ the length of the vector $T^n (v)$ is bounded by
$$ Q(1+\delta )^n |\lambda |^n , $$
where the factor $(1+\delta )^n$ is required to handle additional 
polynomial growth in the presence of a possible
non-trivial generalized eigenspace. In view of the recursion defining 
$M_{m,n}$ we then see that any scalar multiple of $v$ which lies in 
$M_{m,n}$ must be bounded in length by 
$$ Qrm (1+\delta )^{n-1}|\lambda |^{n-1} + Qrm (1+\delta )^{n-2} 
|\lambda |^{n-2} + \cdots + Qrm , $$
which in turn is bounded by
$$ Qrmn (1+\delta )^n \max (|\lambda |^n , 1) . $$
It follows that any scalar multiple of $v$ contained in
$M_{m,0} + M_{m,1} + \cdots + M_{m,n-1}$
is bounded in length by
$$ Qrm \sum_{j=0}^{n-1} j (1+\delta )^j \max (|\lambda |^j , 1) , $$
and this expression is less than
$$ Qrmn^2 (1+\delta )^n \max (|\lambda |^n , 1) . $$
Since the set $M_{m,0} + M_{m,1} + \cdots + M_{m,n-1}$ is a rectangular 
box squarely positioned with respect to the basis 
$\mathcal{B}$ and centred at the origin (as
described above), we combine these length estimates to conclude that 
$$ \vol (M_{m,0} + M_{m,1} + \cdots 
+ M_{m,n-1}) \leq (Qrmn^2)^p (1+\delta )^n 
\prod_{|\lambda_i | \geq 1}|\lambda_i |^n , $$
which yields the result.
\end{proof}

\begin{proposition}\label{P-upperbound}
Suppose $u\in A_\rho$ is a unitary with $L(u) < \infty$. Then
$$ \Entp_L ( \Ad u \circ\alpha_T \circ\gamma_t ) 
\leq\sum_{|\lambda_i | \geq 1} \log |\lambda_i | $$
where $\lambda_1 , \cdots , \lambda_p$
are the eigenvalues of $T$ counted with spectral multiplicity.
\end{proposition}
  
\begin{proof}
Set $\alpha = \Ad u \circ\alpha_T \circ\gamma_t$ for notational brevity.
Let $\Omega\in Pf(\mathcal{L}\cup (A_\rho )_1 )$ and $\delta > 0$. 
By Lemma~\ref{L-convrate} we can find an $C>0$ such that
$$ \| a - \sigma_n (a) \| \leq C \frac{\log n}{n} $$
for all $n\in\mathbb{N}$ and $a\in\Omega\cup\{ u \}$, where $\sigma_n (a)$
is the $n$th Ces\`{a}ro mean, as defined in the paragraph preceding the 
statement of Lemma~\ref{L-integrals}. Since
$\sigma_n (u^* ) = \sigma_n (u)^*$ we also then have
$$ \| u^* - \sigma_n (u^* ) \| \leq C \frac{\log n}{n} $$
for all $n\in\mathbb{N}$. Furthermore 
$$ \| \alpha^j (a) - \alpha^j (\sigma_n (a)) \| \leq C 
\frac{\log n}{n} $$
for all $j,n\in\mathbb{N}$.
By applying the triangle inequality $n$ times in the usual way to estimate
differences of products and using the fact that the operation of
taking a Ces\`{a}ro is norm-decreasing (as can be seen from the first
display in the statement of Lemma~\ref{L-integrals}), we then have, for any 
$a_1 , \dots , a_n \in\Omega$,
$$ \| a_1 \alpha (a_2 ) \cdots \alpha^{n-1}(a_n ) -
\sigma_{n^2}(a_1 ) \alpha (\sigma_{n^2}(a_2 )) \cdots 
\alpha^{n-1}(\sigma_{n^2}(a_n )) \| \leq C \frac{\log n^2}{n} , $$
and this last quantity is less than $\delta$ for all $n$ greater than
or equal to some $n_0 \in\mathbb{N}$.

With the notation of the statement of Lemma~\ref{L-card} we next note that
for any $a\in A$ and $n\in\mathbb{N}$ we have by definition
$$ \sigma_{n^2}(a) \in \text{span}\big\{ u_1^{k_1} u_2^{k_2} 
\cdots u_p^{k_p} : (k_1 , \dots , k_p ) \in K_{n^2} \big\} , $$
while if 
$$ a\in \text{span}\big\{ u_1^{k_1} u_2^{k_2} 
\cdots u_p^{k_p} : (k_1 , \dots , k_p ) \in K \big\} $$
for some finite $K\subset\mathbb{Z}^p$ then
$$ (\Ad u )(\sigma_{n^2}(a)) \in \text{span}\big\{ u_1^{k_1} u_2^{k_2} 
\cdots u_p^{k_p} : (k_1 , \dots , k_p ) \in K + K_{2n^2} \big\} $$
for all $n\in\mathbb{N}$ (the factor of $2$ in the subscript of $K_{2n^2}$
is required to handle multiplication of $a$ by both $u$ and $u^*$). 
Thus, since $\gamma_t$ commutes with the operation of taking a Ces\`{a}ro
sum of a given order, the set of all products $\sigma_{n^2}(a_1 ) \alpha 
(\sigma_{n^2}(a_2 )) \cdots \alpha^{n-1}(\sigma_{n^2}(a_n ))$
with $a_i \in \Omega$ for $i=1, \dots n$ is contained in the subspace
$$ X_n = \text{span}\big\{ u_1^{k_1} u_2^{k_2} 
\cdots u_p^{k_p} : (k_1 , \dots , k_p ) \in L_{2n^2 , 0} +
L_{2n^2 , 1} + \cdots + L_{2n^2 , n-1} \big\} , $$
again using the notation in the statement of Lemma~\ref{L-card} (taking
$m=2n^2$ here). In view of the first paragraph $X_n$ approximately contains
$\Omega\cdot\alpha (\Omega ) \cdot\cdots\cdot\alpha^{n-1}(\Omega )$
to within $2\delta$ for all $n\geq n_0$, and by Lemma~\ref{L-card}
there exists a $Q>0$ such that
$$ \dim (X_n ) \leq (2Qn^3)^p (1+\delta )^n 
\prod_{|\lambda_i | \geq 1} |\lambda_i |^n $$
for all $n\in\mathbb{N}$.   
Therefore
\begin{align*}
\Entp_L (\alpha , \Omega , 2\delta ) &\leq
\limsup_{n\to\infty} \frac1n \log \Bigg((2Qn^3)^p (1+\delta )^n 
\prod_{|\lambda_i | \geq 1} |\lambda_i |^n \Bigg) \\
&= \log (1+\delta ) + \sum_{|\lambda_i | \geq 1} \log |\lambda_i | .
\end{align*} 
Taking the supremum over all $\delta > 0$ then yields
$$ \Entp_L (\alpha , \Omega ) \leq \sum_{|\lambda_i | \geq 1} 
\log |\lambda_i | , $$
from which the proposition follows.
\end{proof}

\begin{theorem}\label{T-toralauto}
We have
$$ \Entp_L (\alpha_T \circ\gamma_t ) = \Entp_{L, \tau}
(\alpha_T \circ\gamma_t  ) =  \sum_{|\lambda_i | \geq 1} 
\log |\lambda_i | $$
where $\lambda_1 , \cdots , \lambda_p$
are the eigenvalues of $T$ counted with spectral multiplicity.
In particular, 
$$ \Entp_L (\alpha_T ) = \Entp_{L, \tau}
(\alpha_T ) = \sum_{|\lambda_i | \geq 1} \log |\lambda_i | . $$
\end{theorem}

\begin{proof}
This follows by combining Propositions \ref{P-lowerbound}, 
\ref{P-upperbound}, and \ref{P-comparison}.
\end{proof}

We also have the following, which is a consequence of 
Propositions \ref{P-comparison} and \ref{P-upperbound}.

\begin{proposition}\label{P-innerauto}
If $u\in A$ is a unitary with $L(u) < \infty$ then 
$$ \Entp_L (\Ad u ) = \Entp_{L, \tau}(\Ad u ) = 0 . $$
\end{proposition}

It is readily seen that if $u$ is a unitary of the form 
$\eta u_1^{k_1} u_2^{k_2}\cdots u_p^{k_p}$ for some integers 
$k_1 , \dots k_p$ and complex number $\eta$ 
of unit modulus, then the automorphism $\Ad u \circ\alpha_T \circ\gamma_t$
can be alternatively expressed as $\alpha_T \circ\gamma_{t'}$ for some 
$t' \in\mathbb{T}^p$, in which case Theorem~\ref{T-toralauto} applies.
We leave open the problem of computing the product entropies 
of $\Ad u \circ\alpha_T \circ\gamma_t$ when $u\in\mathcal{L}$ is a unitary
not of this form and the eigenvalues of $T$ do not all lie on
the unit circle. We expect however that the entropies are positive when 
$\alpha_T$ is asymptotically Abelian (see \cite{EP} for a description
of when this occurs in the case $p=2$) and the partial 
Fourier sums or Ces\`{a}ro means of $u$ converge sufficiently fast 
to $u$, for we could then aim to apply the argument of the 
proof of Proposition~\ref{P-lowerbound} up to a degree of 
approximation.


\begin{thebibliography}{999}

\bibitem{AF}
Alicki, R., and Fannes, M.: Defining quantum dynamical entropy.
{\rm Lett.\ Math.\ Phys.} {\bf 32} (1994), 75--82.

\bibitem{AFTA}
Andries, J., Fannes, M., Tuyls, P., and Alicki, R.: The dynamical entropy
of the quantum Arnold cat map. {\rm Lett.\ Math.\ Phys.} {\bf 35} (1995),
375--383.

\bibitem{Bre} 
Brenken, B.: Representations and automorphisms of the irrational rotation
algebra. {\rm Pacific J. Math.} {\bf 111} (1984), 257--282.

\bibitem{Br} 
Brown, N. P.: Topological entropy in exact $C^*$-algebras.
{\rm Math.\ Ann.} {\bf 314}, 347--367 (1999)

\bibitem{C1}
Connes, A.: Compact metric spaces, Fredholm modules and hyperfiniteness.
{\rm Ergod.\ Th.\ Dynam.\ Sys.} {\bf 9}, 207--220 (1989)

\bibitem{C2}
Connes, A.: {\it Noncommutative Geometry}. San Diego: Academic Press, 1994

\bibitem{C3}
Connes, A.: Gravity coupled with matter and the foundation of
non-commutative geometry. {\rm Commun.\ Math.\ Phys.} {\bf 182}, 155--176
(1996)

\bibitem{CNT} 
Connes, A., Narnhofer, H., and Thirring, W.: Dynamical
entropy of $C^*$-algebras and von Neumann algebras. 
{\rm Commun.\ Math.\ Phys.} {\bf 112} (1987), 691--719.

\bibitem{CS}
Connes, A., and St{\o}rmer, E.: Entropy of automorphisms of II$_1$-von
Neumann algebras. {\rm Acta.\ Math.} {\bf 134} (1975), 289--306.

\bibitem{DGS}
Denker, M., Grillenberger, C., and Sigmund, K.: {\it Ergodic Theory on
Compact Spaces}. Lecture Notes in Math, vol.\ 527. 
Berlin: Springer-Verlag, 1976

\bibitem{Ell}
Elliott, G. A.: On the $K$-theory of the $C^*$-algebra generated 
by a projective
representation of a torsion-free discrete abelian group. In: {\it Operator
Algebras and Group Representations, Vol.\ I}, pp.\ 159--164. Boston: Pitman, 
1984

\bibitem{Ell86}
Elliott, G. A.: The diffeomorphism group of the irrational rotation
$C^*$-algebra. {\rm C. R. Math.\ Rep.\ Acad.\ Sci.\ Canada} {\bf 8},
329--334 (1986)

\bibitem{QTE}
Hudetz, T.: Quantum topological entropy: first steps of a ``pedestrian''
approach. In: {\it Quantum probability \& related topics}, 
pp.\ 237--261. River Edge, NJ: World Scientific, 1993.

\bibitem{Hu}
Hudetz, T.: Topological entropy for appropriately approximated
$C^*$-algebras. {\rm J. Math.\ Phys,} {\bf 35} (1994), 4303--4333.

\bibitem{Kat}
Katznelson, Y.: {\it An Introduction to Harmonic Analysis}, Second 
Edition. New York: Dover Publications, 1976

\bibitem{KL}
Klimek, S. and Le\'{s}niewski, A.: Quantized chaotic dynamics and
non-commutative KS entropy. {\rm Ann.\ Physics} {\bf 248} (1996),
173--198.

\bibitem{KT}
Kolmogorov, A. N., and Tihomirov, V. M.: $\varepsilon$-entropy and
$\varepsilon$-capacity of sets in functional analysis. {\rm Amer.\
Math.\ Soc.\ Trans.\ (2)} {\bf 17}, 277--364 (1961)

\bibitem{SP}
Makarov, B. M., Goluzina, M. G., Lodkin, A. A., and Podkorytov, A. N.:
{\it Selected Problems in Real Analysis}. Translations of Mathematical
Monographs, Vol.\ 107. Providence: AMS, 1992

\bibitem{EP}
Narnhofer, H.: Ergodic properties of automorphisms on the rotation
algebra. {\rm Rep.\ Math.\ Phys.} {\bf 39} (1997), 387--406.

\bibitem{NT}
Narnhofer, H., and Thirring, W.: $C^*$-dynamical systems that are 
asymptotically highly anticommutative. {\rm Lett.\ Math.\ Phys.} {\bf 35} 
(1995), 145--154.

\bibitem{OPT}
Olesen, D., Pedersen, G. K., and Takesaki, M.: Ergodic actions of compact
Abelian groups. {\rm J. Operator Theory} {\bf 3}, 237--269 (1980)

\bibitem{Pes}
Pesin, Ya.\ B.: {\it Dimension Theory in Dynamical Systems: Contemporary 
Views and Applications}. Chicago: The University of Chicago Press, 1997

\bibitem{Pet}
Peters, J.: Entropy on discrete abelian groups. {\rm Adv.\ Math.}
{\bf 33}, 1--13 (1979) 

\bibitem{Rie}
Rieffel, M. A.: Noncommutative tori---a case study of noncommutative
differentiable manifolds. {\it Contemporary Math.} {\bf 105} (1990), 
191--211.

\bibitem{MSACG}
Rieffel, M. A.: Metrics on states from actions of compact groups.
{\rm Doc.\ Math.} {\bf 3}, 215--229 (1998)

\bibitem{MSS}
Rieffel, M. A.: Metrics on state spaces. {\rm Doc.\ Math.} {\bf 4}, 
559--600 (1999)

\bibitem{GHDQMS}
Rieffel, M. A.: Gromov-Hausdorff distance for quantum metric spaces.
arXiv:math.OA/0011063 v2 (2001)

\bibitem{RD} 
Russo, B., and Dye, H. A.: A note on unitary operators
in $C^*$-algebras. {\rm Duke Math. J.} {\bf 33}, 413--416 (1966)

\bibitem{ST} 
Sauvageot, J.-L., and Thouvenot, P.: Une nouvelle
d\'{e}finition de l'entropie dynamique des syst\`{e}mes non-commutatifs.
{\rm Commun.\ Math.\ Phys.} {\bf 145} (1992), 411--423.

\bibitem{Sla}
Slawny, J.: On factor representations and the $C^*$-algebra of canonical 
commutation relations. {\rm Commun.\ Math.\ Phys.} 
{\bf 24}, 151--170 (1972)

\bibitem{NE}
St{\o}rmer, E.: A survey of noncommutative dynamical entropy. 
In: {\it Classification of Nuclear $C^*$-algebras. Entropy in Operator
Algebras}, pp.\ 147--198. Berlin: Springer, 2002.

\bibitem{Th}
Thomsen, K.: Topological entropy for endomorphisms of local $C^*$-algebras.
{\rm Commun.\ Math.\ Phys.} {\bf 164} (1994), 181--193.

\bibitem{Voi} 
Voiculescu, D. V.: Dynamical approximation entropies
and topological entropy in operator algebras. {\rm Commun.\ Math.\ Phys.}
{\bf 170}, 249--281 (1995)

\bibitem{Wat}
Watatani, Y.: Toral automorphisms on irrational rotation algebras.
{\it Math.\ Japon.} {\bf 26} (1981), 479--484.

\bibitem{LADVNA}
Weaver, N.: Lipschitz algebras and derivations of von Neumann algebras.
{\rm J. Funct.\ Anal.} {\bf 139}, 261--300 (1996)

\bibitem{LANT}
Weaver, N.: $\alpha$-Lipschitz algebras on the noncommutative torus. 
{\rm J. Operator Theory} {\bf 39}, 123--138 (1998)

\bibitem{LA}
Weaver, N.: {\it Lipschitz Algebras}. River Edge, NJ: World Scientific, 
1999

\end{thebibliography}
\end{document}